\documentclass{article}

\usepackage{
 amsmath,
 amsxtra,
 amsthm,
 amssymb,
 etex,
 mathrsfs,
 stmaryrd,
pdflscape
 }
\usepackage[all]{xy}
\usepackage{hyperref}

\newtheorem{theorem}{Theorem}[section]
\newtheorem{lemma}[theorem]{Lemma}

\newtheorem{proposition}[theorem]{Proposition}
\newtheorem{corollary}[theorem]{Corollary}

\theoremstyle{definition}
\newtheorem{defn}[theorem]{Definition}
\newtheorem{remark}[theorem]{Remark}

\newcommand{\bd}{\begin{defn}}
\newcommand{\ed}{\end{defn}}
\newcommand{\bl}{\begin{lemma}}
\newcommand{\el}{\end{lemma}}
\newcommand{\bp}{\begin{proposition}}
\newcommand{\ep}{\end{proposition}}
\newcommand{\bt}{\begin{theorem}}
\newcommand{\et}{\end{theorem}}
\newcommand{\bc}{\begin{corollary}}
\newcommand{\ec}{\end{corollary}}
\newcommand{\br}{\begin{remark}}
\newcommand{\er}{\end{remark}}
\newcommand{\ba}{\begin{array}}
\newcommand{\ea}{\end{array}}
\newcommand{\bpf}{\begin{proof}}
\newcommand{\epf}{\end{proof}}

\newcommand{\Z}{\mathbb{Z}}
\newcommand{\Q}{\mathbb{Q}}
\newcommand{\Zp}{\mathbb{Z}_{p}}
\newcommand{\Qp}{\mathbb{Q}_{p}}

\newcommand{\Op}{\mathcal{O}}

\newcommand{\al}{\alpha}
\newcommand{\be}{\beta}
\newcommand{\Ga}{\Gamma}

\newcommand{\la}{\lambda}

\DeclareMathOperator{\Sel}{Sel} \DeclareMathOperator{\Gal}{Gal}
\DeclareMathOperator{\Hom}{Hom} \DeclareMathOperator{\rank}{rank}
\DeclareMathOperator{\corank}{corank}
\DeclareMathOperator{\Ext}{Ext} \DeclareMathOperator{\Tor}{Tor}

\newcommand{\cts}{\mathrm{cts}}

\newcommand{\cyc}{\mathrm{cyc}}
\newcommand{\ch}{\mathrm{char}}
\newcommand{\p}{\mathfrak{p}}

\newcommand{\ot}{\otimes}
\newcommand{\ilim}{\displaystyle \mathop{\varinjlim}\limits}
\newcommand{\plim}{\displaystyle \mathop{\varprojlim}\limits}

\newcommand{\coker}{\mathrm{coker}\,}

\newcommand{\lra}{\longrightarrow}

\newcommand{\ps}[1]{\llbracket #1 \rrbracket}

  \DeclareFontFamily{U}{wncy}{}
  \DeclareFontShape{U}{wncy}{m}{n}{<->wncyr10}{}
  \DeclareSymbolFont{mcy}{U}{wncy}{m}{n}
  \DeclareMathSymbol{\sha}{\mathord}{mcy}{"58}

\begin{document}
\title{Algebraic functional equation for big Galois representations over multiple $\Zp$-extensions}
 \author{
Zeping Hao\footnote{Hua Loo-Keng Center for Mathematical Science, Academy of Mathematics and Systems Science, Chinese Academy of Sciences, Beijing, 100190, P.R.China.
 E-mail: \texttt{zepinghao@amss.ac.cn}} \quad
   Meng Fai Lim\footnote{School of Mathematics and Statistics, Key Laboratory of Nonlinear Analysis and Applications (Ministry of Education),
Central China Normal University, Wuhan, 430079, P.R.China.
 E-mail: \texttt{limmf@ccnu.edu.cn}}}

\date{}
\maketitle

\begin{abstract} \footnotesize
\noindent
We present a general approach to establish algebraic functional equations for big Galois representations over multiple $\Zp$-extensions. Our result is formulated in both Selmer group and Selmer complex settings, and encompasses a broad range of Iwasawa-theoretic scenarios. In particular, our result applies to the triple product of Hida families in both balanced and unbalanced cases, as well as the half-ordinary Rankin-Selberg universal deformations recently studied by the first named author and Loeffler. Our result also significantly generalizes many previously known cases of algebraic functional equations and answers a question of Greenberg.

\medskip
\noindent Keywords and Phrases:  Algebraic functional equation, big Galois representation, multiple $\Zp$-extension.

\smallskip
\noindent Mathematics Subject Classification 2020: 11G05, 11R23, 11S25.
\end{abstract}

\section{Introduction}

The main conjecture in Iwasawa theory posits a precise relationship between a Selmer group and a conjectural $p$-adic $L$-function. This conjectured $p$-adic $L$-function is expected to satisfy a functional equation. Taking into consideration both the main conjecture and the conjectured functional equation, we expect an algebraic relationship between the Selmer groups corresponding to a Galois representation and its Tate dual, which should mirror the functional equation of the $p$-adic $L$-functions. Greenberg was the first to make an in-depth study of this relationship, naming it the ``algebraic functional equation'' (see \cite{G89}). Moreover, he succeeded in establishing this relationship unconditionally, without assuming the main conjecture nor the existence of the $p$-adic $L$-function. We should note that this result of Greenberg only considered Galois representations with coefficients in a finite extension of $\Qp$.

Following this foundational work, there have been attempts extending this result to big Galois representations, i.e., representations with coefficients in a power series ring (for instance, see \cite{JM, JP}). However, these results have often been restricted to specific types of Galois representations, typically those arising from Hida families, or have been established over the cyclotomic $\mathbb{Z}_p$-extension. A general, unifying framework for big Galois representations over multiple $\mathbb{Z}_p$-extensions does not seem to be present in current existing literature. Additionally, it should be noted that the work of \cite{JM, JP} contains certain flaws in their argument, and our current work will address and rectify some of these inaccuracies.

The goal of this paper is to provide this unifying framework.  We prove an algebraic functional equation for big Galois representations, which satisfy a general axiomatic setting (see conditions \textbf{(C1)-(C4)} and \textbf{(R1)}-\textbf{(R2)} below), over multiple $\mathbb{Z}_p$-extensions that contain the cyclotomic one. This approach allows us to treat diverse arithmetic examples$-$ranging from Hida families to products of Hida families and half ordinary Rankin-Selberg deformations (see Section 7)$-$within a single, comprehensive theory. Consequently, our main theorem not only recovers the prior results mentioned above but significantly generalizes them.

%The goal of this paper is to present a unified framework to prove an algebraic functional equation which is applicable not only to individual Galois representations but to families of them (the so-called big Galois representations) over a multiple $\Zp$-extension. A typical example of such a big Galois representation comes from Hida family \cite{Hida86}. We should mention that the algebraic functional equation for a Hida family has previously been established by Jha-Pal \cite{JP} and Jha-Majumdar \cite{JM}, although their works primarily focused on the cases over a cyclotomic $\Zp$-extension or a $\Zp^2$-extension of an imaginary quadratic field. Our result not only covers these prior studies but also further generalized them. Additionally, our results also covers the context of product of Hida families and the half-ordinary Rankin-Selberg universal deformation of Hao-Loeffler \cite{HL}. In particular, in the situation of a triple product of Hida families, our result applies to both the ``balanced'' and ``unbalanced'' cases in the sense of \cite{Hsi, HsiY}.

We now briefly state our main result. Throughout the paper, $p$ will always denote a fixed odd prime. Let $F$ be a number field. We then let $R$ denote the power series ring $\Op\ps{W_1,...,W_m}$, where $\Op$ is the ring of integers of some finite extension of $\Qp$. We consider the datum $\big(T, \{T_v\}_{v|p}\big)_{R,F}$ which satisfy all of the following properties.

\begin{enumerate}
 \item[(\textbf{C1})] $T$ is a
free $R$-module of rank $d$ with a
continuous, $R$-linear $\Gal(\bar{F}/F)$-action which is
unramified outside a finite set of primes of $F$.

 \item[(\textbf{C2})] For each prime $v$ of $F$ above $p$, there is a
distinguished $R[\Gal(\bar{F}_v/F_v)]$-submodule $T_v$ of $T$ such that $T_v$ is a free $R$-module  of $R$-rank $d_v$, and $T_v^- := T/T_v$ is a free $R$-module of
$R$-rank $d-d_v$.

 \item[(\textbf{C3})] For each real prime $v$ of $F$, $T_v^+:=
T^{\Gal(\bar{F}_v/F_v)}$  is free of
$R$-rank $d^+_v$.

\item[(\textbf{C4})] The following equality
  \[\sum_{v|p} (d-d_v)[F_v:\Qp] = dr_2(F) +
 \sum_{v~\mathrm{real}}(d-d^+_v)\]
holds. Here $r_2(F)$ denotes the number of complex primes of $F$.
\end{enumerate}

For a $R[\Gal(\overline{F}/F)]$-module $N$, we write $N(1) = N\ot_{\Zp}\Zp(1)$, where $\Zp(1)$ is the Tate module of the group of all $p$-power roots of unity, and $\Gal(\overline{F}/F)$ acts diagonally on the tensor product. Throughout the paper, we also use $(-)^\vee$ to denote the Pontryagin dual $\Hom_{\mathrm{cts}}(-,\Qp/\Zp)$. With these notations, we can now define the Tate dual datum $\big(T^*, \{T_v^*\}_{v|p}\big)_{R,F}$ of $\big(T, \{T_v\}_{v|p}\big)_{R,F}$, where we set $T^*= \Hom_R(T,R(1))$ and $T_v^* = \Hom_R(T_v^-, R(1))$ for each $v|p$.

Setting $A= T\ot_R R^\vee$ and $A_v= T_v\ot_R R^\vee$, we then call $(A, \{A_v\}_{v|p})_{R, F}$ the discrete datum of $\big(T, \{T_v\}_{v|p}\big)_{R,F}$. The discrete datum $(A^*, \{A_v^*\}_{v|p})_{R, F}$ of $\big(T^*, \{T^*_v\}_{v|p}\big)_{R,F}$ is defined similarly.

Let $F_\infty$ be a $\Zp^r$-extension of $F$ which is always assumed to contain the cyclotomic $\Zp$-extension $F_\cyc$ of $F$. Following Greenberg, we define the Selmer groups for the discrete data over $F_\infty$, thereby obtaining two Selmer groups, denoted as $\Sel_{Gr}(A/F_\infty)$ and $\Sel_{Gr}(A^*/F_\infty)$. We will use $X_{Gr}(A/F_\infty)$ and $X_{Gr}(A^*/F_\infty)$ to denote the Pontryagin dual of the respective Selmer groups. These groups are endowed with module structures over $R\ps{G}$, where $G=\Gal(F_\infty/F)$. Note that $R\ps{G}$ is now isomorphic to a power series ring (over $\Op$) in $(m+r)$-variables. Before presenting our main results, we need to introduce two additional assumptions on our datum (and its dual).

\textbf{(R1)}: The modules $H^0(F_\cyc, A)$ and $H^0(F_\cyc, A^*)$ are cotorsion over $R$.

\textbf{(R2)}: The modules $H^0(F_{\cyc,w}, A_v^-)$ and $H^0(F_{\cyc,w}, (A_v^*)^-)$ are cotorsion over $R$ for every $v$ dividing $p$ and every prime $w$ of $F_\cyc$ above $v$.

Our first main theorem is as follows.

\bt[Theorem \ref{main thm0}]
Suppose that $R=\Op$ and that our datum  $(T, \{T_v\}_{v|p})_{\Op,F}$ satisfies \textbf{(C1)-(C4)}, \textbf{(R1)} and \textbf{(R2)}. Write $\Ga=\Gal(F_\cyc/F)$. Then $X_{Gr}(A/F_\cyc)$ and $X_{Gr}(A^*/F_\cyc)$ have the same $\Op\ps{\Ga}$-ranks, and there is a pseudo-isomorphism
\[X_{Gr}(A/F_\cyc)_{\mathrm{tor}}\sim \big(X_{Gr}(A^*/F_\cyc)_{\mathrm{tor}}\big)^\iota\]
of $\Op\ps{\Ga}$-modules.
\et

In particular, our theorem  provides an affirmative answer to Greenberg's question \cite[P130, (66)]{G89} whether the modules $X_{Gr}(A/F_\cyc)_{\mathrm{tor}}$ and $\big(X_{Gr}(A^*/F_\cyc)_{\mathrm{tor}}\big)^\iota$ are pseudo-isomorphic to each other. We should mention that Lee \cite[Theorem 5.3.3]{Lee} has previously established such a result for abelian varieties with good ordinary reduction at all primes above $p$, via a different approach. In contrast, our argument presented here follows Greenberg's original idea.

Our next main result is concerned with the higher analog of the algebraic functional equation, namely, either when the Galois representation is defined over $R=\Op\ps{W_1,...,W_m}$ with $m\geq 1$ or the $\Zp^r$-extension considered has $r\geq 2$. As it turns out, we can package these two extensions into one single theorem.

\bt (Theorem \ref{main thm})
Let $F_\infty$ be a $\Zp^r$-extension of $F$ which contains the cyclotomic $\Zp$-extension. Write $G=\Gal(F_\infty/F)$. Suppose that our datum $(T, \{T_v\}_{v|p})_{R,F}$ satisfies \textbf{(C1)-(C4)}, as well as \textbf{(R1)} and \textbf{(R2)}. Then the following assertions hold true.

\begin{enumerate}
  \item[$(a)$] $X_{Gr}(A/F_\infty)$ is a torsion $R\ps{G}$-module if and only if $X_{Gr}(A^*/F_\infty)$ is a torsion $R\ps{G}$-module.
\item[$(b)$] $X_{Gr}(A/F_\infty)$ is finitely generated over $R\ps{H}$ if and only if $X_{Gr}(A^*/F_\infty)$ is finitely generated over $R\ps{H}$, where $H=\Gal(F_\infty/F_\cyc)$.
  \item[$(c)$] In the event that aforementioned finite generation assertion holds for at least one (and hence both) of the Selmer groups, we then have the following equality
\[ \ch_{R\ps{G}}\Big(X_{Gr}(A/F_\infty)\Big) = \ch_{R\ps{G}}\Big(X_{Gr}(A^*/F_\infty)^\iota\Big).\]
\end{enumerate}
\et
Here, for a $R\ps{G}$-module $M$, we denote by $M^\iota$ the $R\ps{G}$-module which is $M$ as a $R$-module but with $g\in G$ acting by $g^{-1}$.

 The aforementioned theorems can also be reformulated using the Selmer complexes of Nekov\'a\v{r} \cite{Nek}.
 Writing $SC(T/F_\infty)$ and $SC(T^*/F_\infty)$ for the Selmer complexes corresponding to our given data and its dual, respectively (see Section \ref{Selmer section} for their precise definitions), with cohomology groups $H^i(SC(T/F_\infty))$ and $H^i(SC(T^*/F_\infty))$, Theorems \ref{main thm0} and \ref{main thm}, in the language of Selmer complexes, will take the following forms.

\bt [Theorem \ref{main thm0 Selcomplex}]
Suppose that $R=\Op$ and that our datum $(T, \{T_v\})_{\Op, F}$ satisfies \textbf{(C1)-(C4)}, \textbf{(R1)} and \textbf{(R2)}. Then $H^2(SC(T/F_\cyc))$ and $H^2(SC(T^*/F_\cyc))$ have the same $\Op\ps{\Ga}$-ranks. Moreover, we have a pseudo-isomorphism
\[H^2(SC(T/F_\cyc))_{\mathrm{tor}}\sim \big(H^2(SC(T^*/F_\cyc))_{\mathrm{tor}}\big)^\iota.\]
\et

\bt [Theorem \ref{main thm Selcomplex}]
Let $F_\infty$ be a $\Zp^r$-extension of $F$ which contains the cyclotomic $\Zp$-extension. Suppose that our datum $(T, \{T_v\}_{v|p})_{R,F}$ satisfies \textbf{(C1)-(C4)}, \textbf{(R1)} and \textbf{(R2)}. Then the following statements are valid.
\begin{enumerate}
  \item[$(a)$]  $H^2(SC(T/F_\infty))$ is a torsion $R\ps{G}$-module if and only if $H^2(SC(T^*/F_\infty))$ is a torsion $R\ps{G}$-module.
       \item[$(b)$]  $H^2(SC(T/F_\infty))$ is finitely generated over $R\ps{H}$ if and only if $H^2(SC(T^*/F_\infty))$ is finitely generated over $R\ps{H}$, where $H=\Gal(F_\infty/F_\cyc)$.
  \item[$(c)$] Suppose that the aforementioned finite generation assertion holds for at least one (and consequently, both) of the cohomology groups. Then we have the following:
\[\ch_{R\ps{G}}\Big(H^2(SC(T/F_\infty))\Big) = \ch_{R\ps{G}}\Big( H^2(SC(T^*/F_\infty))^\iota\Big).\]
\end{enumerate}
\et

We make a few remarks concerning the formulation of our main result. Conditions \textbf{(C1)-(C4)} are modelled after the so-called Panchishkin condition (for instance, see \cite[Section 4]{G94}), although our setup only retains the algebraic aspect and deliberately omits other arithmetic conditions, including but not limited to Hodge-Tate properties, purity, and $p$-criticality. The conditions \textbf{(R1)} and \textbf{(R2)} seems to be satisfied for most arithmetic objects of interest (see \cite{CSW, KT}). In Section \ref{examples section}, we will provide a selection of examples (by no means an exhaustive list) to illustrate the applicability of our theorem. The astute readers would have observed that the algebraic functional equation in Theorems \ref{main thm} and \ref{main thm Selcomplex} are obtained under a stringent condition than just merely requiring the Iwasawa modules in question to be torsion. The reason behind is because we make use an algebraic result (see Proposition \ref{alg lemma}) that imposes us such a condition to work with. In \cite{JM, JP}, the authors have stated an algebraic result in this spirit which is stronger than ours. Unfortunately, their results are incorrect (see Remarks \ref{Jha et al} for explanations and counterexamples to their results). Our work not only resolves this particular issue but also generalizes the algebraic equation to  more general big Galois representations.

%In Section \ref{examples section}, we will present a few examples (by no means an exhaustive list) where our theorem applies. It is interesting to note that our result, in particular, applies to abelian variety with good ordinary reduction at all primes above $p$, and $p$-ordinary modular forms. Algebraic functional equations in these context have been established in \cite{JM, JP, Lee} over the cyclotomic $\Zp$-extension. Our theorem thus not only recover these prior results but also generalizes them to multiple $\Zp$-extensions. For big Galois representations, we give examples coming from product of Hida families and the half-ordinary Rankin-Selberg universal deformation of Hao-Loeffler \cite{HL}. In particular, in the situation of a triple product of Hida families, our result can cover both the ``balanced'' and ``unbalanced'' cases. For more in-depth discussion on this, we refer readers to Section \ref{examples section}.

%Another important reason for presenting our result in such a general form is that it greatly facilitates our eventual proof of the algebraic functional equation for big Galois representations.
Without further delay, we shall now describe the outline of the proof leaving details to the body of the paper. The first step of the proof involves introducing an alternative Selmer group defined by Greenberg, known as the strict Selmer group (see \cite[Page 121]{G89}). The motivation for utilizing this group stems from the fact that, while the Greenberg Selmer group is indeed the appropriate object for formulating the main conjecture, the strict Selmer group is more amenable to global duality arguments which suits our purposes for proving Theorem \ref{main thm0}. Furthermore, the strict Selmer groups serve as an intermediary object, enabling us to establish a connection between the Greenberg Selmer group and the Selmer complex (see Remark \ref{strict Sel suffices}).  As it turns out, under our assumption \textbf{(R2)}, we can show that the Pontryagin dual of the Greenberg Selmer group, the Pontryagin dual of the strict Selmer group and the second cohomology group of the Selmer complex are pseudo-isomorphic to one another (see Lemmas \ref{2 Selmer groups} and \ref{strict selmer Selcomplex}). Hence Theorem \ref{main thm0} is equivalent to
Theorem \ref{main thm0 Selcomplex}, and the same applies to Theorem \ref{main thm} and Theorem \ref{main thm Selcomplex}. Consequently, in either context, proving one will be sufficient.

For the proof of Theorem \ref{main thm0} (which is equivalent to Theorem \ref{main thm0 Selcomplex}), we will prove the assertions for the strict Selmer group following Greenberg's approach. For the proof of Theorem \ref{main thm} (or its equivalent, Theorem \ref{main thm Selcomplex}), we prove the Selmer complex version, and the proof proceeds in the following sequence of steps.
For the first step, we leverage the properties of the Selmer complex to show that, in order to prove the theorem for a given datum over a $\Zp^r$-extension, it is sufficient to establish the theorem for another larger datum that is constructed from the original datum but over the cyclotomic $\Zp$-extension. Although this larger datum is defined over a power series ring with more variables than the original, it still falls under our axiomatic setting. The crux of this reduction lies in the fact that our proof can now proceed within the fixed setting of $F_\infty =F_\cyc$, thus paving the way for a proof by induction, namely, on the number of variables in the power series ring $R$ (see Remark \ref{reduction remark} and the preceding discussion).

For the final inductive argument, the idea is to choose infinite many linear ideals (see Definition \ref{linear element definition}) of $\Op\ps{W_1,..., W_m}$ appropriately so that the datum, a priori defined over $\Op\ps{W_1,...,W_m}$, specializes to data defined over $\Op\ps{W_1,...,W_{m-1}}$, whose Selmer complex's second cohomology group is pseudo-isomorphic to the corresponding specialization of the second cohomology group of Selmer complex of the starting datum. It is precisely at this juncture that the generality of our framework proves advantageous, as such specializations may not inherently have any arithmetic meanings, even when the original datum we started with comes from an arithmetic object. Nevertheless, we can show that these specialized data, upon an appropriate choice of linear ideals,  still satisfy \textbf{(C1)-(C4)}, as well as \textbf{(R1)} and \textbf{(R2)}, thus allowing us to conclude that Theorem \ref{main thm} is valid for these specialized data by our induction hypothesis. Subsequently, we make use of an algebraic result (see Proposition \ref{alg lemma}) to obtain the conclusion for the original datum which is defined over $\Op\ps{W_1,..., W_m}$.

We conclude this introductory section by giving an outline of the paper. In Section \ref{Algebraic section}, we review certain algebraic results that are required for our purposes. Moving on to Section \ref{Selmer section}, this is where we introduce the Greenberg Selmer groups along with the Selmer complexes, and it is also in this section that we present our main results. In Section \ref{compare Selmer}, we introduce the strict Selmer group of Greenberg, and use it to relate the Greenberg Selmer group with the Selmer complex. Sections \ref{m=0 section} and \ref{Selmer proof} will be devoted to the proof of  \ref{main thm0}/\ref{main thm0 Selcomplex} and Theorem \ref{main thm}/\ref{main thm Selcomplex} respectively. Following that, we then discuss some classes of examples in Section \ref{examples section}. Finally, Section \ref{appendix} supplies certain auxiliary results on the irreducibility of representations coming from product of modular forms that are required for discussing some of the examples in Section \ref{examples section}.

\subsection*{Acknowledgement}

The authors would like to thank Yangyu Fan, Somnath Jha, Antonio Lei, David Loeffler, Chao Qin, Xin Wan, Jun Wang and Sarah Zerbes for their interest and comments on the paper. M. F. Lim's research is partially supported by the Fundamental Research Funds for the Central Universities No. CCNU25JCPT031.

\section{Algebraic preliminaries} \label{Algebraic section}

We collect certain algebraic results that will be required for the discussion of the paper.
Firstly, we make a general observation.
\bl \label{coprime torsion corollary}
Let $R$ be an integral domain, and $\p$ a prime ideal of $R$. Suppose that $M$ is a finitely generated $R$-module such that $M/\p$ is torsion over $R/\p$. Then the following assertions are valid.
 \begin{enumerate}
   \item[$(i)$] $M$ is torsion over $R$.
   \item[$(ii)$] Furthermore, in the event that $R$ is Noetherian, we have that $\Tor_j^R(R/\p, M)$ is torsion over $R/\p$ for every $j\geq 1$.
 \end{enumerate}

\el

\bpf
 We start by claiming that $\p + \mathrm{Ann}_R(M) \neq \p$, and we will set out to establish our claim. Suppose on the
contrary that $\p + \mathrm{Ann}_R(M) = \p$. For an ideal $J$ of $R$, we write $\mathrm{rad}(J)$ for the (prime) radical of $J$, i.e., $\mathrm{rad}(J)$ consists of elements $x\in R$ such that $x^{n(x)}\in J$ for some positive integer $n(x)$. Then we have
\begin{equation}\label{AnnRad} \mathrm{rad}\big(\mathrm{Ann}_{R}(M/\p)\big) = \mathrm{rad}\big(\p + \mathrm{Ann}_{R}(M)\big) =\mathrm{rad} (\p) = \p,\end{equation}
 where the first equality follows from \cite[P.\ 13, Exercise 2.2]{Matsu}, the second by our supposition and the final equality follows from the fact that $\p$ is a prime ideal. On the other hand, in view of the hypothesis that $M/\p$ is torsion over $R/\p$, there exists an element $r\notin \p$ such that $(r+\p)\in R/\p$ annihilates $M/\p$. This in turn implies that $r$ annihilates $M/\p$. However, by (\ref{AnnRad}), the latter will force $r\in \p$ which yields the required contradiction. This completes the proof of the claim.

 By the claim, it particularly follows that $\mathrm{Ann}_R(M)\neq 0$, which thereby establishes (i). For (ii), we first note that $\Tor_R^j(R/\p, M)$ can be viewed as a $R$-module and a $R/\p$-module at the same time. Since $R$ is Noetherian and $M$ is finitely generated over $R$, $\Tor_R^j(R/\p, M)$ is finitely generated over $R/\p$. In view of the latter finite generation, we are reduced to showing that $\mathrm{Ann}_{R/\p}\big(\Tor^R_j(R/\p, M)\big)\neq 0$. Now, viewing $\Tor_R^j(R/\p, M)$ as a $R$-module, we have the following inclusion
 \[ \p + \mathrm{Ann}_R(M) \subseteq \mathrm{Ann}_R\big(\Tor^R_j(R/\p, M)\big). \]
 Combining this observation with our claim, we see that there exists an element $r\in \mathrm{Ann}_R\big(\Tor^R_j(R/\p, M)\big)$ with $r\notin \p$. This in turn implies that $\mathrm{Ann}_{R/\p}\big(\Tor^R_j(R/\p, M)\big)\neq 0$, thereby proving assertion (ii) of the lemma. \epf

The next lemma is essentially taken from \cite{LimFine}.
\bl \label{coprime torsion}
Let $x$ be a prime element of a regular local ring $R$ such that $R/x$ is also a regular local ring. For every finitely generated $R$-module $M$, we then have
\[ \rank_{R}(M) = \rank_{R/x}(M/x) -\rank_{R/x}(M[x]).\]
   In particular, the following two statements are equivalent.
 \begin{enumerate}
  \item[$(a)$] The module $M/x$ is torsion over $R/x$.
   \item[$(b)$] The module $M[x]$ is torsion over $R/x$ and the module $M$ is torsion over $R$.
 \end{enumerate}
 Moreover, if $M$ is torsion over $R$ and $x$ is a prime element of $R$ which does not lie in the support of $M$. Then $M[x]$ and $M/x$ are torsion modules over the ring $R/x$.
\el

\bpf
The equality of ranks follows from \cite[Proposition 4.12]{LimFine}. The equivalence of the two statements is then immediate from this. For the final assertion, a commutative algebraic argument tells us that $M/x$ is a torsion $R/x$-module (for instance, see \cite[Lemma 2.2]{LimFineMod}). By the asserted equality of ranks of the proposition, we then see that $M[x]$ is also torsion over $R/x$.
\epf

For the remainder of this section, we shall let $R$ denote the power series ring $\Op\ps{W_1,...,W_m}$ in $m$-variables with coefficients in $\Op$, where $\Op$ is the ring of integers of some fixed finite extension of $\Qp$. We also fix a uniformizer $\pi$ of $\Op$.
The following definitions are taken from Ochiai \cite{Och}.

\bd (\cite[Definition 3.2]{Och}) \label{linear element definition}
\begin{itemize}
  \item[(1)] An element $\ell$ of $R$ is said to be a linear element if it is of the form $a_0+a_1W_1+\cdots +a_mW_m$ with $a_i\in\Op$, where $a_0$ is divisible by $\pi$ and not all $a_i$ are divisible by $\pi$. An ideal $I$ of $R$ is said to be a linear ideal if there exists a linear element $\ell$ such that $I=(\ell)$. Note that a linear ideal is necessarily a prime ideal of $R$ of height $1$.
        %\item[(2)] An automorphism $\sigma$ of $R$ is said to be a linear transform if $\sigma(W_j) = \sum_{i=1}^mt_{i,j}W_i$ with $(t_{i,j})_{1\leq i,j\leq m}\in \mathrm{GL}_m(\Op)$.
  \item[(2)] Let $m\geq 2$. For a finitely generated torsion $R$-module $M$, we denote by $\mathcal{L}_R(M)$ the collection of all linear ideals $(\ell)$ such that $M/\ell$ is a torsion $R/\ell$-module and that the image of $\ch_{R}(M)$ in $R/\ell$ is equal to $\ch_{R/\ell}(M/\ell)$.
\end{itemize}
\ed

\bl \label{linear element lemma}
\begin{itemize}
  \item[$(i)$] For every linear element $\ell$, the quotient ring $\Op\ps{W_1,...,W_m}/\ell$ is isomorphic to a power series ring in $(m-1)$-variables with coefficients in $\Op$.

  %\item[$(ii)$] Let $\sigma$ be a linear transform, and $x\in R$. Then $x$ is a linear element if and only if $\sigma(x)$ is a linear element.
  \item[$(ii)$] Let $\ell,\ell'$ be two linear elements of $R$. Then $(\ell)= (\ell')$ if and only if $\ell = u\ell'$ for some $u\in \Op^\times$.
\end{itemize}
\el

\bpf
For a given linear element $\ell$, upon relabeling, we may assume $a_m$ is not divisible by $\pi$, and so the map
      \[ \Op\ps{W_1,...,W_m}\lra \Op\ps{W_1,...,W_{m-1}},\quad W_m\mapsto - a_m^{-1}(a_0+a_1W_1+\cdots +a_{m-1}W_{m-1})\]
      will induces an isomorphism $\Op\ps{W_1,...,W_m}/\ell\cong \Op\ps{W_1,...,W_{m-1}}$.
      This proves statement (i). Statement (ii) is \cite[Lemma 3.4]{Och}.
\epf

We are now in position to state the following result which will be required for the eventual proof of our main theorem. %This is inspired by the result of Ochiai \cite[Proposition 3.6]{Och}, also see \cite[Proposition 4.7]{JP} or \cite[Proposition 3.6]{JM}, but see the below remark.
In what follows, we denote by $\mathfrak{m}$ the maximal ideal of $R=\Op\ps{W_1,...,W_m}$. Let $W$ be another free variable. A $R\ps{W}$-module $M$ may naturally be viewed as a $R$-module. Furthermore, we have $\mathfrak{m}M = (\mathfrak{m}R\ps{G})M$, where $\mathfrak{m}R\ps{G}$ is the ideal of $R\ps{G}$ generated by $\mathfrak{m}$. Therefore, we write $M/\mathfrak{m}$ for $M/\mathfrak{m}M = M/\mathfrak{m}R\ps{G}M$, which can be viewed as a $\Op/\pi\ps{W}$-module. Also, we will often view a linear element $\ell$ (resp., linear ideal $(\ell_i)$) of $R$ as a linear element (resp., linear ideal) of $R\ps{W}$ without further mention.

The following lemma is standard, and we record it here for the convenience of the readers.

\bl \label{torsion psuedo-null}
Let $M$ be a $R\ps{W}$-module which is finitely generated as a $R$-module. Then $M$ is pseudo-null over $R\ps{W}$ if and only $M$ is torsion over $R$.
\el

\bpf
See \cite[Lemma 5.1]{LimFine}.
\epf

We now prove the following important algebraic result which will be a key ingredient to the proof of our main result.

\bp \label{alg lemma}
Let $M$ and $N$ be finitely generated torsion $R\ps{W}$-modules, where $R=\Op\ps{W_1,...,W_m}$ with $m\geq 1$.  Suppose that the following statements are valid.
\begin{enumerate}
\item[$(a)$] The modules $M$ and $N$ are finitely generated as $R$-modules.
\item[$(b)$] There exist infinitely many distinct linear ideals $(\ell_i)$ of $R$ satisfying all of the following.
\begin{enumerate}
    \item[$(i)$] For every $i$, we have $(\ell_i)\in \mathcal{L}_{R\ps{W}}(M)\cap \mathcal{L}_{R\ps{W}}(N)$.
   \item[$(ii)$] For every $i$, there is an identity $\ch_{R_i\ps{W}}(M/\ell_i) = \ch_{R_i\ps{W}}(N/\ell_i)$, where $R_i=R/\ell_i$.
 \end{enumerate}
  \end{enumerate}
Then we have $\ch_{R\ps{W}}(M) = \ch_{R\ps{W}}(N)$.
\ep

\bpf
 %We begin proving $\rank_{R}(M) = \rank_{R}(N)$, or equivalently, $a=b$. Choose an $\ell_0$ that satisfies (a)-(c) and does not lie in the support of $M$ and $N$. Then by \cite[Lemma 2.2]{LimFineMod}, we have that $M_{R-\mathrm{tor}}/\ell_0$ and $N_{R-\mathrm{tor}}/\ell_0$ are torsion over $R_i$. By Lemma \ref{coprime torsion}, this in turn implies that $M_{R-\mathrm{tor}}[\ell_0]$ and $N_{R-\mathrm{tor}}[\ell_0]$ are torsion over $R_0:=R/\ell_0$. But since one has $M[\ell_0] = M_{R-\mathrm{tor}}[\ell_0]$ and $N[\ell_0] = N_{R-\mathrm{tor}}[\ell_0]$, we see that $M[\ell_0]$ and $N[\ell_0]$ are torsion over $R_0$. Then by appealing to Lemma \ref{coprime torsion} again, we obtain $\rank_R(M) = \rank_{R_0}(M/\ell_0)$ and $\rank_R(N) = \rank_{R_0}(N/\ell_0)$, which upon combining with (b), yields the required conclusion.
 Set $f: = \ch_{R\ps{W}}(M)$ and $g: = \ch_{R\ps{W}}(N)$. Write $f = \sum_{n\geq 0}c_n W^n$, where $c_n\in R$. We first show that not all of the $c_n$'s lie in $\mathfrak{m}$. Given that $M$ is torsion over $R\ps{W}$, there is a short exact sequence
 \[ 0\lra R\ps{W}/f_1\times \cdots \times R\ps{W}/f_s \lra M \lra P \lra 0, \]
 with $f= \prod_{j=1}^s f_j$ and $P$ being a pseudo-null $R\ps{W}$-module. From which, we obtain an exact sequence
 \[ \Tor_{R\ps{W}}^1(R\ps{W}/\mathfrak{m}, P) \lra R\ps{W}/(f_1, \mathfrak{m})\times \cdots \times R\ps{W}/(f_s, \mathfrak{m})\lra M/\mathfrak{m} \lra P/\mathfrak{m} \lra 0\]
 As $M$ is finitely generated over $R$ by assumption (a), it follows from Nakayama lemma that $M/\mathfrak{m}$ is finite which in turn implies that $P/\mathfrak{m}$ is finite. By Lemma \ref{coprime torsion corollary}, we then have that $\Tor_{R\ps{W}}^1(R\ps{W}/\mathfrak{m}, P) $ is also finite. Hence we may conclude that $R\ps{W}/(f_j, \mathfrak{m})$ is finite for every $j$. But $R\ps{W}/(f_j, \mathfrak{m})\cong \Op/\pi\ps{W}/ \overline{f}_j$, where $\overline{f}_j = f_j$ (mod $\mathfrak{m}$), and the latter is finite if and only if $\overline{f}_j\neq 0$, or equivalently, $f_j\notin \mathfrak{m}R\ps{W}$. Since $\mathfrak{m}R\ps{W}$ is a prime ideal of $R\ps{W}$, we see that $f\notin \mathfrak{m}R\ps{W}$, and this proves what we set to show.

Taking the preceding assertion into account, it then follows from an application of the division lemma (see \cite[5.3.1]{NSW}) that $R\ps{W}/f$ is a free $R$-module of finite rank. On the other hand, since $R$ is a unique factorization domain, there is a natural injection
\begin{equation}\label{UFD}
R/l_1\cdots l_k \hookrightarrow R/\ell_1\times \cdots \times R/\ell_k
\end{equation}
for every $k$. As $R\ps{W}/f$ is a free $R$-module of finite rank, upon applying $\ot_R R\ps{W}/f$, we obtain an injection
\begin{equation}\label{flat tensor}
R\ps{W}/(f, l_1\cdots l_k) \hookrightarrow R\ps{W}/(f, l_1)\times \cdots \times R\ps{W}/(f, l_k)
\end{equation}
Combining this with assumption (b), we see that the image of $g$ in $R\ps{W}/(f, l_1\cdots l_k)$ is zero for every $k$. Letting $k\to\infty$, we have that the image of $g$ in $R\ps{W}/(f)$ is zero.
A similar argument then also shows that the image of $f$ in $R\ps{W}/(g)$ is zero. This completes the proof of our result.
\epf

We now give some remarks concerning the result of Ochiai \cite[Proposition 3.6]{Och} (also see \cite[Proposition 4.7]{JP} or \cite[Proposition 3.6]{JM}).

\br \label{Jha et al}
\begin{enumerate}
\item[$(1)$] Assumption (a) in Proposition \ref{alg lemma} is necessary for the conclusion to hold. We shall give an example to illustrate this. Let $M= \Zp\ps{W_1, W}/(W_1-p)\times \Zp\ps{W_1, W}/(W_1-p)$ and $N= \Zp\ps{W_1,W}/(W_1-p^2)$. Plainly, $\ch_{\Zp\ps{W_1, W}}(M) = (W_1-p)^2 \neq W_1-p^2= \ch_{\Zp\ps{W_1, W}}(N)$. But if we take $\ell_i=W_1-p^{i+1}$ for $i\geq 1$, then one sees easily that $M/\ell_i = \Zp\ps{W}/(p^{i+1}-p)\times \Zp\ps{W}/(p^{i+1}-p)$ and $N/\ell_i = \Zp\ps{W}/(p^{i+1}-p^2)$ with \[\ch_{\Zp\ps{W}}(M/\ell_i) = (p^{i+1}-p)^2 = (p)^2 = (p^{i+1}-p^2)= \ch_{\Zp\ps{W}}(N/\ell_i)\] for every $i\geq 1$. In particular, this also shows that \cite[Proposition 4.7]{JP} and \cite[Proposition 3.6]{JM} are not correct as stated.

    We shall say a bit on the error in the proof of \cite[Proposition 4.7]{JP} and \cite[Proposition 3.6]{JM}. A crucial step in our proof is that the characteristic element $f$ has the property that $R\ps{W}/f$ is finite flat over $R$ (which forces us to impose assumption (a)). In the previous works \cite{JM, JP}, the authors overcome this by applying a change of coordinates, say $(W_1,...,W_m,W)\mapsto (W'_1,...,W'_m,W')$, where the characteristic element $f$ has the property that $\Op\ps{W'_1,...,W'_m,W'}/f$ is finite flat over $\Op\ps{W'_1,...,W'_m}$. However, they fail to take into account that the change of coordinates also affect the linear elements. In particular, the changes of coordinates may cause the linear elements to be not defined over $\Op\ps{W'_1,...,W'_m}$, and so we may not have an analogous injection over $\Op\ps{W'_1,...,W'_m}$ as in (\ref{UFD}) which is necessary for the next line of argument in (\ref{flat tensor}).

\item[$(2)$] In the proof of \cite[Proposition 3.6]{Och} (also see \cite[Proposition 4.7]{JP} or \cite[Proposition 3.6]{JM}), there is a line which says ``the element $\pi^{\mu}$ is equal to the highest power of $\pi$ dividing the characteristic power series of $M/\ell$.'' This line of reasoning is not correct as can be seen in the following example. Indeed, if we write $M=\Zp\ps{W_1,W}/(W_1-p)$ and $\ell = W_1-p^i$ for $i\geq 1$. Then $\ch_{\Zp\ps{W_1, W}}(M) = W_1-p$. Under the projection map $\Zp\ps{W_1, W} \twoheadrightarrow \Zp\ps{W_1, W}/\ell_i$, this ideal is sent to $(p^i-p) = (p)$ which is precisely $\ch_{\Zp\ps{W_1, W}/\ell_i}(M/\ell_i)$ contrary to the aforementioned line.
\end{enumerate}
\er

Finally, in view of assumption (b)(i) of Proposition \ref{alg lemma}, we conclude the section with the following lemma which will be useful for our subsequent discussion.

\bl \label{linear element choice}
Let $M$ be a torsion $R\ps{W}$-module, where $R=\Op\ps{W_1,...,W_m}$ with $m\geq 1$. Denote by $M_{\mathrm{null}}$ the maximal pseudo-null $R\ps{W}$-submodule of $M$. Suppose that $\ell$ is a linear element of $R$ such that $M_{\mathrm{null}}[\ell]$ is pseudo-null over $R/\ell\ps{G}$ and that $M/\ell$ is torsion over $R/\ell\ps{W}$. Then $\ell \in \mathcal{L}_{R\ps{W}}(M)$.
\el

\bpf
For simplicity, we shall write $\overline{R}=R/\ell$. Since $M/\ell$ is already assumed to be torsion over $\overline{R}\ps{W}$, it remains to show that the image of $\ch_{R\ps{W}}(M)$ in $\overline{R}\ps{W}$ is equal to $\ch_{\overline{R}\ps{W}}(M/\ell)$. Consider the following exact sequence
\[ 0 \lra M_{\mathrm{null}}\lra M \lra E \lra N\lra 0, \]
where $N$ is some pseudo-null $R\ps{W}$-module and $E = \bigoplus_{\mathfrak{p}}R\ps{W}/\mathfrak{p}^{a_\mathfrak{p}}$ for some height one prime $\mathfrak{p}$ of $R\ps{W}$ not equal to $(\ell)$ (this is necessary due to the $\overline{R}\ps{W}$-torsionness of $M/\ell$). Denote by $Z$ the image of the map $M\lra E$, we then have two short exact sequences
\[ 0\lra M_{\mathrm{null}}\lra M  \lra Z\lra 0, \]
\[ 0\lra Z\lra E \lra N\lra 0.\]
Since the prime ideals appearing in $E$ are not equal to $\ell$ and $R\ps{W}$ is a unique factorization domain, we have $E[\ell]=0$ which in turn implies that $Z[\ell]=0$. Taking these into account, we obtain the following two exact sequences
\[ 0\lra M_{\mathrm{null}}/\ell\lra M/\ell  \lra Z/\ell\lra 0, \]
\[ 0\lra N[\ell] \lra Z/\ell\lra E/\ell \lra N/\ell \lra 0.\]
Now since $N$ is a pseudo-null $R\ps{W}$-module, it follows from \cite[Lemma 3.1]{Och} that $\ch_{\overline{R}\ps{W}}(N[\ell])= \ch_{\overline{R}\ps{W}}(N/\ell)$. Thus, it follows from this and the above exact sequence that we have $\ch_{\overline{R}\ps{W}}(Z/\ell])= \ch_{\overline{R}\ps{W}}(E/\ell)$. On the other hand, since $M_{\mathrm{null}}[\ell]$ is pseudo-null over $\overline{R}\ps{W}$, \cite[Lemma 3.1]{Och} also tells us that $M_{\mathrm{null}}/\ell$ is pseudo-null over $\overline{R}\ps{W}$. Therefore, the above short exact sequence yields
$\ch_{\overline{R}\ps{W}}(M/\ell)= \ch_{\overline{R}\ps{W}}(Z/\ell)$, and hence we have the identity
\[ \ch_{\overline{R}}(M/\ell)= \ch_{\overline{R}\ps{W}}(E/\ell) = \ch_{R\ps{W}}(E)/\ell.\]
This is precisely what we set out to show.
\epf

\section{Main results} \label{Selmer section}

Let $p$ be an odd prime, and $F$ a number field. For the remainder of the paper, $R$ will always denote $\Op\ps{W_1,...,W_m}$  which is a power series ring in $m$-variables (we allow $m=0$) with coefficient in $\Op$, where $\Op$ is the ring of integers of some fixed finite extension of $\Qp$.

Even though this has already been introduced in the introductory section, for the readers' convenience, we will once again restate the axiomatic framework pertaining to our datum, which, as before, will be denoted by $\big(T, \{T_v\}_{v|p}\big)_{R,F}$.

\begin{enumerate}
 \item[(\textbf{C1})] $T$ is a
free $R$-module of rank $d$ with a
continuous, $R$-linear $\Gal(\bar{F}/F)$-action which is
unramified outside a finite set of primes of $F$.

 \item[(\textbf{C2})] For each prime $v$ of $F$ above $p$, there is a
distinguished $R[\Gal(\bar{F}_v/F_v)]$-submodule $T_v$ of $T$ such that $T_v$ is a free $R$-module  of
$R$-rank $d_v$, and $T_v^- := T/T_v$ is a free $R$-module of
$R$-rank $d-d_v$.

 \item[(\textbf{C3})] For each real prime $v$ of $F$, $T_v^+:=
T^{\Gal(\bar{F}_v/F_v)}$  is free of
$R$-rank $d^+_v$.

\item[(\textbf{C4})] The following equality
  \[\sum_{v|p} (d-d_v)[F_v:\Qp] = dr_2(F) +
 \sum_{v~\mathrm{real}}(d-d^+_v)\]
holds. Here $r_2(F)$ denotes the number of complex primes of $F$.
\end{enumerate}

As we need to work with Selmer groups defined over a tower of number fields, we need to consider the base change of our datum. For a finite extension $L$ of $F$, the base change of the datum $\big(T,
\{T_w\}_{w|p} \big)_{R,L}$ over $L$ is given as follows:

(1) We view $T$ as a $\Gal(\overline{F}/L)$-module via restriction of the Galois action.

(2) For each prime $w$
of $L$ above $p$, we set $T_w =T_v$, where $v$ is the prime of $F$
below $w$, and view it as a $\Gal(\bar{F}_v/L_w)$-module via the appropriate restriction.
Note that we then have $d_w= d_v$.

(3) For each real prime $w$ of $L$ which lies above a real prime $v$ of $F$, we set
$T_w^+= T^{\Gal(\bar{F}_v/F_v)}$ and write $d^+_w
= d^+_v$.

The following lemma gives some
sufficient conditions for the equality in \textbf{(C4)} to hold for
the datum $\big(T, \{T_w\}_{w|p}\big)_{R,L}$.

\bl \label{data base change} Let $\big(T, \{T_v\}_{v|p}\big)_{R,F}$ be a datum defined over $F$, which satisfies \textbf{(C1)-(C4)}. Suppose that $L$ is a finite extension of $F$ and that at least one of the following statements holds:
\begin{enumerate}
  \item[$(i)$] All the archimedean primes of $F$ are unramified in $L$
  \item[$(ii)$] $L$ is a Galois extension of $F$ with $[L:F]$ being odd.
  \item[$(iii)$] $F$ is totally imaginary.
  \item[$(iv)$] $F$ is totally real, $L$ is totally imaginary and
 \[ \sum_{v~\mathrm{real}}d_v^+ = d[F:\Q]/2.\]
\end{enumerate}
 Then the datum $\big(T, \{T_w\}_{w|p}\big)_{R,L}$ obtained via base changing to $L$ also satisfies \textbf{(C1)-(C4)}. \el

\bpf It suffices to check that the datum $\big(T, \{T_w\}_{w|p}\big)$ satisfies \textbf{(C4)}. As a start, we note that the validity of (ii) or (iii) implies the validity of (i). Thus, to establish the lemma in these three cases, it suffices to prove it assuming the validity of (i). For this, we have the following calculation
\[ \ba{rl}
\displaystyle \sum_{w|p} (d-d_w)[L_w:\Qp]\!\!
   &= \displaystyle\sum_{v|p}\sum_{w|v} (d-d_v)[L_w:F_v][F_v:\Qp] \\
   &= \displaystyle\sum_{v|p} (d-d_v)[F_v:\Qp] \sum_{w|v}[L_w:F_v] \\
   &= \displaystyle [L:F]\sum_{v|p} (d-d_v)[F_v:\Qp] \\
  &= \displaystyle d[L:F]r_2(F) + [L:F]\sum_{v~\mathrm{real}} (d-d^+_v). \\  \ea \]
By virtue of assumption (i), every real prime (resp., complex prime) of $F$ remains a real prime (resp., complex prime) in $L$. Therefore, we have $[L:F] r_2(F) = r_2(L)$ and
\[ [L:F]\sum_{v~\mathrm{real}} (d-d^+_v) = \sum_{w~\mathrm{real}} (d-d^+_w). \]
This proves what we want to show.

Now suppose that (iv) holds. Then we have $r_2(F)=0$ and
\[ \ba{rl}
\displaystyle\sum_{w|p} (d-d_w)[L_w:\Qp]\!\!
& = \displaystyle\sum_{v|p}\sum_{w|v} (d-d_v)[L_w:F_v][F_v:\Qp] \\
&= \displaystyle\sum_{v|p}(d-d_v)[F_v:\Qp]\sum_{w|v} [L_w:F_v] \\
&= [L:F] \displaystyle\sum_{v~\mathrm{real}}(d-d_v^+)\\
&= [L:F] \displaystyle\sum_{v~\mathrm{real}}d - [L:F] \displaystyle\sum_{v~\mathrm{real}}d_v^+\\
&= [L:F] [F:\Q] d - [L:F] [F:\Q]d/2 \\
&= d[L:\Q]/2 = dr_2(L). \\
\ea \]
Taking into account that $L$ has no real primes, this proves \textbf{(C4)}.
  \epf

For a $R[\Gal(\overline{F}/F)]$-module $N$, we write $N(1) = N\ot_{\Zp}\Zp(1)$, where $\Zp(1)$ is the Tate module of the group of all $p$-power roots of unity, and $\Gal(\overline{F}/F)$ acts diagonally on the tensor product. With this notation, we can now define the Tate dual datum of $\big(T, \{T_v\}_{v|p}\big)_{R,F}$, where we set $T^*= \Hom_R(T,R(1))$, and for each $v|p$ (resp. real prime $v$), we define $T_v^* = \Hom_R(T_v^-, R(1))$ (resp., $(T_v^*)^+ = \Hom_R(T/T^+_v, R(1))$).

Denote by $A$ the discrete Galois representation $T\ot_R\Hom_{\mathrm{cont}}(R,\Qp/\Zp)$. We then write $A_v$ (resp., $A_v^-$) for
$T_v\ot_R\Hom_{\mathrm{cont}}(R,\Qp/\Zp)$ (resp., $T_v^-\ot_R\Hom_{\mathrm{cont}}(R,\Qp/\Zp)$). We will call $(A, \{A_v\}_{v|p})_{R,F}$ the discrete datum attached to $\big(T, \{T_v\}_{v|p}\big)_{R,F}$.
The discrete datum $(A^*, \{A^*_v\}_{v|p})_{R,F}$ for $\big(T^*, \{T^*_v\}_{v|p}\big)_{R,F}$ is defined similarly.

The relationship between $T, T^*, A$ and $A^*$ is illustrated by the following diagram
\[  \entrymodifiers={!! <0pt, 0.8ex>+} \SelectTips{eu}{}\xymatrix@R=.6in @C=1in{
     T \ar[d]_{-\ot_{R}R^\vee} \ar[r]^{\Hom_{R}(-, R(1))} \ar[dr] &
    T^* \ar[l] \ar[dl] \ar[d]^{-\ot_{R}R^\vee}
     \\
    A  \ar[ur]_{}
    &  A^* \ar[ul]}  \]
where the two diagonal maps are given by $\Hom_{\cts}(-,\mu_{p^{\infty}})$.

\br
Throughout the rest of this paper, we shall call $\big(T, \{T_v\}_{v|p}\big)_{R,F}$ and/or $(A, \{A_v\}_{v|p})_{R,F}$ as our ``datum''. We then call $\big(T^*, \{T^*_v\}_{v|p}\big)_{R,F}$ and/or $(A^*, \{A^*_v\}_{v|p})_{R,F}$ as our ``dual datum''. This terminological convention will be adhered to in the subsequent discussion of the paper without further mention. Finally, in contexts where the meaning is unambiguous, we simplify our notation by omitting $R$ and $F$, and simply write $\big(T, \{T_v\}_{v|p}\big)$, $(A, \{A_v\}_{v|p})$ and similar expressions for the dual datum.
\er

 Denoting by $F_\cyc$ the cyclotomic $\Zp$-extension of $F$, we impose two final assumptions on our datum and its dual.

\textbf{(R1)}: The modules $H^0(F_\cyc, A)$ and $H^0(F_\cyc, A^*)$ are cotorsion over $R$.

\textbf{(R2)}: The modules $H^0(F_{\cyc,w}, A_v^-)$ and $H^0(F_{\cyc,w}, (A_v^*)^-)$ are cotorsion over $R$ for every $v$ dividing $p$ and every prime $w$ of $F_\cyc$ above $v$.

\medskip
We now recall the definition of the Greenberg Selmer group. As a start, we let $S$ be a finite set of primes of $F$ that contains all the primes of $F$ above $p$, the ramified primes of $T$ and the infinite primes. The maximal algebraic extension of $F$ unramified outside $S$ is then denoted by $F_S$. For every algebraic extension $\mathcal{F}$ of $F$ contained in $F_S$, we write $G_S(\mathcal{F}) =\Gal(F_S/\mathcal{F})$.
Let $F_{\infty}$ be a $\Zp^r$-extension of $F$ which is always assumed to contain the cyclotomic $\Zp$-extension $F_\cyc$. We allow $r=1$, in which case, we then have $F_\infty=F_\cyc$. Note that one always had $F_\infty\subseteq F_S$ (see \cite[Theorem 1]{Iw73}). Let $L$ be a finite extension of $F$ contained in $F_\infty$. Consider the following local condition
\[ H^1_{Gr}(L_w, A)=
\begin{cases} \ker\big(H^1(L_w, A)\lra H^1(L_w^{ur}, A^-_w)\big) & \text{\mbox{if} $w|p$},\\
 \ker\big(H^1(L_w, A)\lra H^1(L^{ur}_w, A)\big) & \text{\mbox{if} $w\nmid p$,}
\end{cases} \]
where $L_w^{ur}$ is the maximal unramified extension of $L_w$. The Greenberg Selmer group attached to the datum $\big(A,
\{A_w\}_{w|p} \big)_{R,L}$  is then given by
\[ \Sel_{Gr}(A/L) = \ker\left( H^1(G_S(L),A)\lra \bigoplus_{v \in S}\bigoplus_{w|v}\frac{H^1(L_w, A)}{H^1_{Gr}(L_w,
A)}\right).\]
We define $\Sel_{Gr}(A/F_{\infty}) = \ilim_L \Sel_{Gr}(A/L)$,
where the limit runs over all finite extensions $L$ of $F$ contained
in $F_{\infty}$. The Pontryagin dual of $\Sel_{Gr}(A/F_{\infty})$ is then denoted by $X_{Gr}(A/F_{\infty})$, and this group comes equipped with a $R\ps{G}$-module structure, where $G=\Gal(F_{\infty}/F)$.

For a $R\ps{G}$-module $M$, we denote by $M^\iota$ the $R\ps{G}$-module which is $M$ as a $R$-module but with $g\in G$ acting by $g^{-1}$. With this, we can now state the main theorem of this paper.

\bt \label{main thm0}
Suppose that $R=\Op$ and that our datum $(T, \{T_v\})_{\Op, F}$ satisfies \textbf{(C1)-(C4)}, \textbf{(R1)} and \textbf{(R2)}. Then $X_{Gr}(A/F_\cyc)$ and $X_{Gr}(A^*/F_\cyc)$ have the same $\Op\ps{\Ga}$-ranks. Moreover, we have a pseudo-isomorphism
\[X_{Gr}(A/F_\cyc)_{\mathrm{tor}}\sim \big(X_{Gr}(A^*/F_\cyc)_{\mathrm{tor}}\big)^\iota,\]
where we denote $M_{\mathrm{tor}}$ the torsion submodule of a given $R\ps{\Ga}$-module $M$.
\et

We now make the following remark.

\br
In \cite[P130, (66)]{G89}, Greenberg raised the question whether the modules $X_{Gr}(A/F_\cyc)_{\mathrm{tor}}$ and $\big(X_{Gr}(A^*/F_\cyc)_{\mathrm{tor}}\big)^\iota$ are pseudo-isomorphic to each other. Our Theorem \ref{main thm0} thus provides an affirmative answer to Greenberg's question. It is worth noting that Lee \cite[Theorem 5.3.3]{Lee}) previously established this result for abelian varieties with good ordinary reduction at all primes above $p$, via a different approach. In contrast, our argument follows Greenberg's original idea.
\er

%Let $\mathfrak{m}$ denote the maximal ideal of $R$. Consider the residual (discrete) datum $(\widetilde{A},(\widetilde{A}_v)_{v|p})$, where $\widetilde{A}:=A[\mathfrak{m}]$ and $\widetilde{A}_v:= A_v[\mathfrak{m}]$. We then denote by the residual Greenberg Selmer group $\Sel_{Gr}(\widetilde{A}/F_\cyc)$ over $F_\cyc$, whose definition is analogous to that of the Greenberg Selmer group, with $A$ replaced by $\widetilde{A}$.

For big Galois representations over multiple $\Zp$-extensions, our result for their algebraic functional equation is as follows.

\bt \label{main thm}
Let $F_\infty$ be a $\Zp^r$-extension of $F$ which contains the cyclotomic $\Zp$-extension. Suppose that our datum $(T, \{T_v\}_{v|p})_{R,F}$ satisfies \textbf{(C1)-(C4)}, \textbf{(R1)} and \textbf{(R2)}. Then the following statements are valid.

\begin{enumerate}
  \item[$(a)$] $X_{Gr}(A/F_\infty)$ is a torsion $R\ps{G}$-module if and only if $X_{Gr}(A^*/F_\infty)$ is a torsion $R\ps{G}$-module.
      \item[$(b)$] $X_{Gr}(A/F_\infty)$ is finitely generated over $R\ps{H}$ if and only if $X_{Gr}(A^*/F_\infty)$ is finitely generated over $R\ps{H}$, where $H=\Gal(F_\infty/F_\cyc)$.
  \item[$(c)$] In the event that (b) holds, then we have the following equality
\[ \ch_{R\ps{G}}\Big(X_{Gr}(A/F_\infty)\Big) = \ch_{R\ps{G}}\Big(X_{Gr}(A^*/F_\infty)^\iota\Big).\]
\end{enumerate}
\et

We now recast our results in the language of the Selmer complexes as developed by Nekov\'a\v{r} \cite{Nek}. We begin recalling the definition of the Selmer complex for the Greenberg local conditions. For a profinite group $\mathcal{G}$ and a topological $\mathcal{G}$-module $N$, we let $C(\mathcal{G}, N)$ denote the complex of continuous cohains on $\mathcal{G}$ with values in $N$ with the usual differential maps. Following \cite[Subsections 6.7 and 7.8]{Nek}, for each $v\in S$, we define
\[ U_v(T) = \left\{
                  \begin{array}{ll}
                    C(\Gal(\bar{F}_v/F_v), T_v), & \hbox{if $v|p$;} \\
                    C(\Gal(F_v^{ur}/F_v), (T)^{I_v}), & \hbox{if $v\nmid p$,}
                  \end{array}
                \right.
  \]
where $I_v$ is the inertia group of $\Gal(\bar{F}_v/F_v)$. There is a natural chain map $i_v: U_v(T^*) \lra C(\Gal(\bar{F}_v/F_v), T)$.
Write $U_S(T)=\bigoplus_{v\in S} U_v(T)$ and $i_S^+=\bigoplus_{v\in S}i_v$. The Selmer complex attached to the datum $(T, \{T_v\}_{v|p}\}$ is then defined to be
\[ SC(T/F) = \mathrm{cone}\left(C(G_S(F), T)\oplus U_S(T^*) \stackrel{\mathrm{res}-i_S^+}{\lra} \bigoplus_{v\in S} C(\Gal(\bar{F}_v/F_v), T)  \right)[-1].\]
For every finite extension $L$ of $F$ contained in $F_\infty$, we can similarly define $SC(T/L)$. We then set $SC(T/F_\infty)$ to be $\plim_L SC(T/L)$ as in the sense of  Nekov\'a\v{r} \cite[Subsections 8.6-8.8]{Nek}.

Theorems \ref{main thm0} and \ref{main thm}, in the language of Selmer complexes, will take the following forms.

\bt \label{main thm0 Selcomplex}
Suppose that $R=\Op$ and that our datum $(T, \{T_v\})_{\Op, F}$ satisfies \textbf{(C1)-(C4)}, \textbf{(R1)} and \textbf{(R2)}. Then $H^2(SC(T/F_\cyc))$ and $H^2(SC(T^*/F_\cyc))$ have the same $\Op\ps{\Ga}$-ranks. Moreover, we have a pseudo-isomorphism
\[H^2(SC(T/F_\cyc))_{\mathrm{tor}}\sim \big(H^2(SC(T^*/F_\cyc))_{\mathrm{tor}}\big)^\iota.\]
\et

\bt \label{main thm Selcomplex}
Let $F_\infty$ be a $\Zp^r$-extension of $F$ which contains the cyclotomic $\Zp$-extension. Suppose that our datum $(T, \{T_v\}_{v|p})_{R,F}$ satisfies \textbf{(C1)-(C4)}, \textbf{(R1)} and \textbf{(R2)}. Then the following statements are valid.
\begin{enumerate}
  \item[$(a)$]  $H^2(SC(T/F_\infty))$ is a torsion $R\ps{G}$-module if and only if $H^2(SC(T^*/F_\infty))$ is a torsion $R\ps{G}$-module.
       \item[$(b)$]  $H^2(SC(T/F_\infty))$ is finitely generated over $R\ps{H}$ if and only if $H^2(SC(T^*/F_\infty))$ is finitely generated over $R\ps{H}$, where $H=\Gal(F_\infty/F_\cyc)$.
  \item[$(c)$] Suppose that the aforementioned finite generation assertion holds for at least one (and consequently, both) of the cohomology groups. Then we have the following:
\[\ch_{R\ps{G}}\Big(H^2(SC(T/F_\infty))\Big) = \ch_{R\ps{G}}\Big( H^2(SC(T^*/F_\infty))^\iota\Big).\]
\end{enumerate}
\et

The proof of the Theorems \ref{main thm0} and \ref{main thm0 Selcomplex} will be given in Section \ref{m=0 section}, whereas Theorems \ref{main thm} and \ref{main thm Selcomplex} will be proven in Section \ref{Selmer proof}.

\section{Comparing Selmer groups and Selmer complexes} \label{compare Selmer}
In this section, we relate the Greenberg Selmer groups and Selmer complexes. To facilitate this discussion, we introduce another Selmer group of Greenberg, specifically, the strict Selmer group, as defined in \cite[Page 121]{G89}. Let $L$ be a finite extension of $F$ which is contained in $F_\infty$. Now, consider the following local condition
\[ H^1_{str}(L_w, A)=
\begin{cases} \ker\big(H^1(L_w, A)\lra H^1(L_w, A^-_w)\big) & \text{\mbox{if} $w|p$},\\
 \ker\big(H^1(L_w, A)\lra H^1(L^{ur}_w, A)\big) & \text{\mbox{if} $w\nmid p$.}
\end{cases} \]
The strict Selmer group attached to the datum $\big(A,
\{A_w\}_{w|p} \big)_{R,L}$  is then defined by
\[ \Sel_{str}(A/L) = \ker\left( H^1(G_S(L),A)\lra \bigoplus_{v \in S}\bigoplus_{w|v}\frac{H^1(L_w, A)}{H^1_{str}(L_w,
A)}\right).\]
We set $\Sel_{str}(A/F_{\infty}) = \ilim_L \Sel_{str}(A/L)$,
where the limit runs over all finite extensions $L$ of $F$ contained
in $F_{\infty}$. The Pontryagin dual of $\Sel_{str}(A/F_{\infty})$ is then denoted by $X_{str}(A/F_{\infty})$.

\bl \label{localcondition}
Let $F_{\infty}$ be a $\Zp^r$-extension of $F$ which contains $F_\cyc$.
For each $v\in S$, we fix a prime of $F_\infty$ above $v$ which we, by abuse of notation, denote by $v$. Writing $G_v$ for the decomposition group of $G$ at $v$, one then has the following identifications
\[ \ilim_L\bigoplus_{w|v}\frac{H^1(L_w, A)}{H^1_{Gr}(L_w,
A)} \cong\left\{
           \begin{array}{ll}
             \mathrm{Coind}^G_{G_v} \Big(H^1(F_{\infty, v}^{ur}, A^-_v)^{\Gal(F_{\infty,v}^{ur}/F_{\infty,v})}\Big), & \hbox{if $v\mid p$} \\
             \mathrm{Coind}^G_{G_v} \Big(H^1(F_{\infty,v}, A)\Big), & \hbox{if $v\nmid p$,}
           \end{array}
         \right.
  \]
\[ \ilim_L\bigoplus_{w|v}\frac{H^1(L_w, A)}{H^1_{str}(L_w,
A)} \cong\left\{
           \begin{array}{ll}
             \mathrm{Coind}^G_{G_v} \Big(H^1(F_{\infty, v}, A^-_v)\Big), & \hbox{if $v\mid p$} \\
             \mathrm{Coind}^G_{G_v} \Big(H^1(F_{\infty,v}, A)\Big), & \hbox{if $v\nmid p$.}
           \end{array}
         \right.
  \]
\el

\bpf
We shall only supply a proof for the $Gr$ local condition in the case of $v\mid p$, as the remaining cases can be proven similar. By definition, we have the exact sequence
\[  0 \lra H^1_{Gr}(L_v, A)\lra H^1(L_v, A) \lra H^1(L_v^{ur}, A^-_v)\]
which upon taking limit yields
\begin{equation}\label{Gr sesA}
  0 \lra \ilim H^1_{Gr}(L_v, A)\lra H^1(F_{\infty,v}, A) \lra H^1(F_{\infty, v}^{ur}, A^-_v)
\end{equation}
Note that the map $H^1(F_{\infty,v}, A) \lra H^1(F_{\infty, v}^{ur}, A^-_v)$ is given by the composition
\[H^1(F_{\infty,v}, A) \lra H^1(F_{\infty,v}, A^-_v) \lra H^1(F_{\infty, v}^{ur}, A^-_v)^{\Gal(F_{\infty,v}^{ur}/F_{\infty,v})}\subseteq H^1(F_{\infty, v}^{ur}, A^-_v). \]
The first map is surjective by \cite[Theorem 7.1.8]{NSW}, and  since $\Gal(F_{\infty,v}^{ur}/F_{\infty,v})$ has $p$-cohomological dimension at most one, so is the second map. Taking these observations into account, the sequence (\ref{Gr sesA}) in turn induces the following short exact sequence
\begin{equation*}
  0 \lra \ilim H^1_{Gr}(L_v, A)\lra H^1(F_{\infty,v}, A) \lra H^1(F_{\infty, v}^{ur}, A^-_v)^{\Gal(F_{\infty,v}^{ur}/F_{\infty,v})}\lra 0,
\end{equation*}
which is what we want to show.
\epf

For notational brevity, we shall write
\[J_v^{Gr}(A/F_\infty)=\left\{
           \begin{array}{ll}
             \mathrm{Coind}^G_{G_v} \big(H^1(F_{\infty, v}^{ur}, A^-_v)^{\Gal(F_{\infty,v}^{ur}/F_{\infty,v})}\big), & \hbox{if $v\mid p$} \\
             \mathrm{Coind}^G_{G_v} \big(H^1(F_{\infty,v},A)\big), & \hbox{if $v\nmid p$.}
           \end{array}
         \right.\]
 \[J_v^{str}(A/F_\infty)=\left\{
           \begin{array}{ll}
             \mathrm{Coind}^G_{G_v} \big(H^1(F_{\infty, v}, A^-_v)\big), & \hbox{if $v\mid p$} \\
             \mathrm{Coind}^G_{G_v} \big(H^1(F_{\infty,v},A)\big), & \hbox{if $v\nmid p$.}
           \end{array}
         \right.\]

The next lemma gives a relation between the two Selmer groups of Greenberg which will play a crucial role in our subsequent discussion.

\bl \label{2 Selmer groups}
There is a surjection
   \[ X_{Gr}(A/F_\infty)\lra X_{str}(A/F_\infty), \]
   whose kernel is a torsion $R\ps{H}$-module (and hence pseudo-null $R\ps{G}$-module).
  %Therefore, in particular, $X_{str}(A/F_\infty)$ is torsion over $\Op\ps{\Ga}$ if and only if $X_{Gr}(A/F_\infty)$ is too.
  We also have an analogous assertion upon replacing all occurrence of $A$ by $A^*$.
\el

\bpf
By Lemma \ref{torsion psuedo-null}, the pseudo-nullity assertion will follow once we can show the $R\ps{H}$-torsionness property. In the course of the proof, for a given $\Gal(\overline{F}_v /F_v)$-module $Z$, we shall frequently write $Z(\mathcal{L}) = Z^{\Gal(\overline{F}_v /\mathcal{L})}$ for any algebraic extension $\mathcal{L}$ of $F_v$. By definition, the two Selmer groups fit into the following commutative diagram
\[ \entrymodifiers={!! <0pt, .8ex>+} \SelectTips{eu}{}\xymatrix{
     0 \ar[r] &  \Sel_{str}(A/F_\infty)\ar[r] \ar[d]_{} &  H^1(G_S(F_\infty),A) \ar[r] \ar@{=}[d] &  \bigoplus_{v\in S} J_v^{str}(A/F_\infty) \ar[d] \\
     0 \ar[r] &  \Sel_{Gr}(A/F_\infty)\ar[r]_{} &   H^1(G_S(F_\infty),A) \ar[r] &  \bigoplus_{v\in S}J_v^{Gr}(A/F_\infty)
     }\]
with exact rows. Taking Lemma \ref{localcondition} into account, we see that the kernel of the rightmost vertical map is given by the finite sum
\[\bigoplus_{v\mid p}\mathrm{Coind}^{G}_{G_v}H^1\big(\Gal(F_{\infty,v}^{ur}/F_{\infty,v}),  A^-_v(F_{\infty,v}^{ur})\big). \]
By the snake lemma, we then have the following exact sequence
\begin{equation}\label{2 Selmer sequence}
0 \lra \Sel_{str}(A/F_\infty) \lra \Sel_{Gr}(A/F_\infty)
  \lra \bigoplus_{v\mid p}\mathrm{Coind}^{G}_{G_v}H^1\big(\Gal(F_{\infty,v}^{ur}/F_{\infty,v}),  A^-_v(F_{\infty,v}^{ur})\big).
\end{equation}
Now observe that for each $v|p$, we have
\[ \mathrm{Coind}^{G}_{G_v}H^1\big(\Gal(F_{\infty,v}^{ur}/F_{\infty,v}),  A^-_v(F_{\infty,v}^{ur})\big) = \bigoplus_{w|v}\mathrm{Coind}^{H}_{H_w}H^1\big(\Gal(F_{\infty,v}^{ur}/F_{\infty,v}),  A^-_v(F_{\infty,v}^{ur})\big), \]
where $w$ runs through the (finitely many) primes of $F_\cyc$ above $v$.
  It therefore remains to show that each
\[H^1\big(\Gal(F_{\infty,v}^{ur}/F_{\infty,v}),  A^-_v(F_{\infty,v}^{ur})\big)^\vee\]
is torsion over $R\ps{H_w}$. In fact, we're going to show an even stronger result: this module is torsion over $R$. Denote by $F_{\infty,v}^{ur,p}$ the maximal $p$-extension of $F_{\infty, v}$ contained in $F_{\infty,v}^{ur}$. Then we have
  \[H^1\big(\Gal(F_{\infty,v}^{ur}/F_{\infty,v}),  A^-_v(F_{\infty,v}^{ur})\big)\cong H^1\big(\Gal(F_{\infty,v}^{ur,p}/F_{\infty,v}),  A^-_v(F_{\infty,v}^{ur,p})\big).\]
   Now, had $\Gal(F_{\infty,v}^{ur,p}/F_{\infty,v})$ being finite, then $H^1\big(\Gal(F_{\infty,v}^{ur,p}/F_{\infty,v}),  A^-_v(F_{\infty,v}^{ur,p})\big)^\vee$ is annihilated by a sufficiently large power of $p$, and so is necessarily a torsion $R$-module. Thus, we may assume that $\Gal(F_{\infty,v}^{ur,p}/F_{\infty,v})\cong \Zp$. Since $A^-_v(F_{\infty,v}^{ur,p})^\vee$ is finitely generated over $R$, it is in particular a torsion $R\ps{\Gal(F_{\infty,v}^{ur,p}/F_{\infty,v})}$-module. In view of this, Lemma \ref{coprime torsion} then yields the following equality
  \[ \corank_R\Big(H^1(\Gal(F_{\infty,v}^{ur,p}/F_{\infty,v}),  A^-_v(F_{\infty,v}^{ur,p}))\Big) =
   \corank_R\Big(H^0(\Gal(F_{\infty,v}^{ur,p}/F_{\infty,v}),  A^-_v(F_{\infty,v}^{ur,p}))\Big).  \]
  But $H^0(\Gal(F_{\infty,v}^{ur,p}/F_{\infty,v}),  A^-_v(F_{\infty,v}^{ur,p}))$ is precisely $H^0(F_{\infty,w}, A_v^-)$ which is cotorsion over $R$ by our assumption \textbf{(R2)}.
\epf

We now relate the strict Selmer group to the Selmer complex.
To do this, we recall from \cite[Definition 5.3.1.1]{Nek} the complex of continuous
cochains with compact support
\[ C_c(G_S(F), T^*) = \mathrm{cone}\left(C(G_S(F), T^*) \stackrel{\mathrm{res}}{\lra} \bigoplus_{v\in S} C(\Gal(\bar{F}_v/F_v), T^*)  \right)[-1]\]
which fits into the following short exact sequence
\[ 0 \lra C_c(G_S(F), T^*) \lra SC(T^*/F) \lra U_S(T^*) \lra 0 \]
of complexes (see \cite[(6.1.3)]{Nek}). Upon taking cohomology, we obtain the following long exact sequence
\begin{multline} \label{long selmer seq}
\bigoplus_{v|p}H^0(F_v, T_v^*)\oplus \bigoplus_{v\nmid p} H^0(\Gal(F_v^{ur}/F_v), (T^*)^{I_v}) \lra  H^1_c(G_S(F), T^*) \lra H^1(SC(T^*/F)) \\
  \lra\bigoplus_{v|p}H^1(F_v, T_v^*)\oplus \bigoplus_{v\nmid p} H^1(\Gal(F_v^{ur}/F_v), (T^*)^{I_v}) \stackrel{\rho}{\lra} H^2_c(G_S(F), T^*)  \\
 \lra H^2(SC(T^*/F)) \lra \bigoplus_{v|p}H^2(F_v, T_v^*) \lra H^3_c(G_S(F), T^*) \lra H^3(SC(T^*/F)) \lra 0.
\end{multline}
By Poitou-Tate duality and local Tate duality (for instances, see \cite[Propositions 5.2.4 and 5.2.7]{Nek} or \cite[Theorem 4.2.6]{LimPT}), the above exact sequence can be rewritten as
\begin{multline} \label{long selmer seq2}
\bigoplus_{v|p}H^0(F_v, T_v^*)\oplus \bigoplus_{v\nmid p} H^0(\Gal(F_v^{ur}/F_v), (T^*)^{I_v}) \lra  H^2(G_S(F),A)^\vee \lra H^1(SC(T^*/F)) \\
  \lra\bigoplus_{v|p}H^1(F_v, A_v^-)^\vee\oplus \bigoplus_{v\nmid p} H^1(\Gal(F_v^{ur}/F_v), A^{I_v})^\vee \stackrel{\rho}{\lra} H^1(G_S(F), A)^\vee  \\
 \lra H^2(SC(T^*/F)) \lra \bigoplus_{v|p}H^0(F_v, A_v^-)^\vee \lra H^3_c(G_S(F), T^*) \lra H^3(SC(T^*/F)) \lra 0.
\end{multline}
In particular, the Pontryagin dual of the two terms in the middle line of the above, together with their map $\rho$, identifies with
\[  H^1(G_S(F), A) \stackrel{\rho^\vee}{\lra} \bigoplus_{v|p}H^1(F_v, A_v^-)\oplus \bigoplus_{v\nmid p} H^1(\Gal(F_v^{ur}/F_v), A^{I_v}).  \]
This in turn implies that $\coker \rho = X_{str}(A/F)$. Therefore, we obtain the following exact sequence
\begin{multline*}
 0\lra X_{str}(A/F) \\ \lra  H^2(SC(T^*/F)) \lra \bigoplus_{v|p}H^0(F_v,A_v^-)^\vee
 \lra H^0(G_S(F), A)^\vee \lra H^3(SC(T^*/F)) \lra 0.
\end{multline*}
For a finite extension $L$ of $F$ contained in $F_\infty$, we have a similar exact sequence as above, and upon taking inverse limit, we obtain
\begin{multline} \label{long sequence} 0\lra X_{str}(A/F_\infty)  \lra  H^2(SC(T^*/F_\infty)) \lra \bigoplus_{v|p}\mathrm{Coind}^G_{G_v} \Big(H^0(F_{\infty, v}, A^-_v)\Big)^\vee  \\
\lra H^0(G_S(F_\infty), A)^\vee \lra H^3(SC(T^*/F_\infty)) \lra 0,
\end{multline}
where $SC(T^*/F_\infty)$ is the inverse limit of the complexes $SC(T^*/L)$ in the sense of \cite{Nek}. We then have the following lemma.

\bl \label{strict selmer Selcomplex}
Let $F_\infty$ be a $\Zp^r$-extension of $F$ which contains the cyclotomic $\Zp$-extension. Suppose that our datum $(T, \{T_v\}_{v|p})_{R,F}$ satisfies \textbf{(C1)-(C4)}, \textbf{(R1)} and \textbf{(R2)}. Then the following statements are valid.
\begin{enumerate}
  \item[$(a)$]  There is an injection $X_{str}(A/F_\infty)\lra H^2(SC(T^*/F_\infty))$ with cokernel being finitely generated torsion over $R\ps{H}$.
       \item[$(b)$]  Both $H^3(SC(T/F_\infty))$ and $H^3(SC(T^*/F_\infty))$ are torsion over $R\ps{H}$ (and so pseudo-null over $R\ps{G}$).
\end{enumerate}
\el

\bpf
By the long exact sequence (\ref{long sequence}), it suffices to show that $H^0(G_S(F_\infty), A)^\vee$ and
\[ \bigoplus_{v|p}\mathrm{Coind}^G_{G_v} \Big(H^0(F_{\infty, v}, A^-_v)\Big)^\vee \]
are finitely generated torsion $R\ps{H}$-modules. Now, if $F_\infty=F_\cyc$, and hence $H=1$, then $H^0(G_S(F_\cyc), A)^\vee$ is indeed torsion over $R$ by \textbf{(R1)}. In the event that $F_\infty\neq F_\cyc$, and hence $H\cong \Zp^{r-1}$ with $r-1\geq 1$, then $H^0(G_S(F_\infty), A)^\vee$ is at most finitely generated over $R$ and so is necessarily torsion over $R\ps{H}$.
For the local cohomology groups, a similar argument appealing to \textbf{(R2)} applies.
\epf

\br \label{strict Sel suffices}
Taking the preceding lemma and Lemma \ref{2 Selmer groups} into account, we see that Theorem \ref{main thm0} and Theorem \ref{main thm0 Selcomplex} are equivalent. Therefore, proving either of them will suffice. Likewise, the same applies to Theorems \ref{main thm} and \ref{main thm Selcomplex}.

Although we do not require it, we thought it might be of interest to write down the following exact sequence
\begin{multline*}  \bigoplus_{v\mid p}\mathrm{Coind}^{G}_{G_v}H^1\big(\Gal(F_{\infty,v}^{ur}/F_{\infty,v}),  A^-_v(F_{\infty,v}^{ur})\big)^\vee \lra X_{Gr}(A/F_\infty)  \lra  H^2(SC(T^*/F_\infty)) \\ \lra \bigoplus_{v|p}\mathrm{Coind}^G_{G_v} \Big(H^0(F_{\infty, v}, A^-_v)\Big)^\vee
\lra H^0(G_S(F_\infty), A)^\vee \lra H^3(SC(T^*/F_\infty)) \lra 0,
\end{multline*}
which relates the Greenberg Selmer group to the Selmer complex. This said exact sequence is obtained by splicing the dual of the sequence (\ref{2 Selmer sequence}) with (\ref{long sequence}).
 \er

\section{Proof of Theorems \ref{main thm0} and \ref{main thm0 Selcomplex}} \label{m=0 section}

In this section, we shall prove Theorems \ref{main thm0} and \ref{main thm0 Selcomplex}. Here, as before, $\Op$ denotes the ring of integers of a finite extension of $\Qp$. Note that in this context, the conditions \textbf{(R1)} and \textbf{(R2)} translate into the following two assertions.

\textbf{(R1)}: The groups $H^0(G_S(F_\cyc), A)$ and $H^0(G_S(F_\cyc), A^*)$ are finite.

\textbf{(R2)}: The groups $H^0(F_{\cyc,w}, A_v^-)$ and $H^0(F_{\cyc,w}, (A_v^*)^-)$ are finite for every $v$ dividing $p$ and every prime $w$ of $F_\cyc$ above $v$.

In the forthcoming discussion of this section, we will denote $\Ga$ as $\Gal(F_\cyc/F)$ and use $M_{\mathrm{tor}}$ to denote the torsion submodule of $M$, where $M$ is a given $\Op\ps{\Ga}$-module. Additionally, we will occasionally, without explicit mention, identify $\Op\ps{\Ga}$ with $\Op\ps{W}$, a power series ring in one variable. We also fix a uniformizer $\pi$ of $\Op$.

In \cite[Theorem 2]{G89}, Greenberg established an algebraic functional equation for datum coming from ``motive of pure weight'' which is $p$-ordinary and $p$-critical. Our proof will follow the approach of Greenberg. Given the way our datum is set up, it may not have such an arithmetic interpretation, and so we shall give a proof here, ensuring that Greenberg's argument follows.

%\bt \label{main thm0}
%Suppose that our datum $(T, \{T_v\})$ satisfies \textbf{(C1)-(C4)}, \textbf{(R1)} and \textbf{(R2)}. Then $X_{Gr}(A/F_\cyc)$ and $X_{Gr}(A^*/F_\cyc)$ have the same $\Op\ps{\Ga}$-ranks. Moreover, we have a pseudo-isomorphism
%\[X_{Gr}(A/F_\cyc)_{\mathrm{tor}}\sim \big(X_{Gr}(A^*/F_\cyc)_{\mathrm{tor}}\big)^\iota.\]
%\et

As a start, we present the following lemma.

\bl \label{structure compare}
Let $M$ and $N$ be two finitely generated $\Op\ps{W}$-modules. Suppose that the following
two statements hold.
\begin{enumerate}
  \item[$(I)$] For every positive integer $n$, we have $\mu_{\Op\ps{\Ga}}(M/\pi^n) =\mu_{\Op\ps{\Ga}}(N/\pi^n)$.
  \item[$(II)$] For every irreducible Weierstrass polynomial $f\in \Op\ps{W}$, we have
  \[ \corank_{\Op}\left((M^\vee\ot_{\Op}\Op\ps{W}/f^n)^\Ga\right)= \corank_{\Op}\left((N^\vee\ot_{\Op}\Op\ps{W}/f^n)^\Ga\right)\]
  for all  positive integers $n$.
\end{enumerate}
Then we have that $\rank_{\Op\ps{W}}(M) = \rank_{\Op\ps{W}}(N)$, and that $M_{\mathrm{tor}}$ is pseudo-isomorphic to $N_{\mathrm{tor}}$.
\el

\bpf
This is essentially contained in \cite[Section 3]{G89} but not explicitly stated there. For the convenience of the readers, we give a sketch of the proof. By \cite[Proposition 2.4.7]{LimCMu}, the first assumption of the lemma will imply that both $\rank_{\Op\ps{W}}(M)$ and $\rank_{\Op\ps{W}}(N)$ are equal, and that $M[\pi^\infty]$ is pseudo-isomorphic to $N[\pi^\infty]$. It therefore remains to show that  $M[f^\infty]$ is pseudo-isomorphic to $N[f^\infty]$ for every irreducible Weierstrass polynomial $f\in \Op\ps{W}$. Suppose that the elementary representations of $M[f^\infty]$ and $N[f^\infty]$ are given by
\[ \bigoplus_{i\geq 1} (\Op\ps{W}/f^i)^{\oplus a_i} \quad\mbox{and} \quad  \bigoplus_{i\geq 1} (\Op\ps{W}/f^i)^{\oplus b_i} \]
respectively, where $a_i=0$ and $b_i=0$ for $i\gg 0$. Let us choose $\theta$ to be a sufficiently large integer so that $a_i=0 = b_i$  for every $i\geq \theta+1$. By the discussion in \cite[Page 109-110]{G89}, we have
\[ \corank_{\Op} \big((M^\vee\ot_{\Op}\Op\ps{W}/f^n)^\Ga\big) =  \corank_{\Op} \Hom_{\Op\ps{W}}(M, K/\Op\ot_{\Op}\Op\ps{W}/f^n).\]
The latter quantity identifies with the $\Op$-rank of $\Hom_{\Op\ps{W}}(M, \Op\ps{W}/f^n)$. In view of this, along with the fact that the functor $\Hom_{\Op\ps{W}}(-, \Op\ps{W}/f^n)$ preserves finite direct sums and sends finite modules to finite modules, one can directly verify that
\[r_n:=\rank_{\Op}\big(\Hom_{\Op\ps{W}}(M, \Op\ps{W}/f^n)\big)= \Big(n \rank_{\Op\ps{W}}(M) + a_1 + 2a_2 +\cdots + (n-1)a_n + n(a_n+\cdots +a_\theta)\Big)\deg(f),\]
\[s_n:=\rank_{\Op}\big(\Hom_{\Op\ps{W}}(N, \Op\ps{W}/f^n)\big)=\Big(n \rank_{\Op\ps{W}}(N) + b_1 + 2b_2 +\cdots + (n-1)b_n + n(b_n+\cdots +b_\theta)\Big)\deg(f)\]
for every $1\leq n\leq \theta$.
Taking into account that we have already shown that $\rank_{\Op\ps{W}}(M)= \rank_{\Op\ps{W}}(N)$, and given that $r_i=s_i$, as stipulated by assumption (II), we can obtain the following matrix equation
\[
 \begin{pmatrix}
  1 & 1 & 1& \cdots &1 \\
  1 & 2 & 2& \cdots &2\\
  1 & 2 & 3& \cdots &3\\
  \vdots  & \vdots & \vdots & \ddots & \vdots\\
  1 & 2 & 3& \cdots & \theta
 \end{pmatrix}
  \begin{pmatrix}
  a_1 \\
  a_2\\
  a_3\\
  \vdots \\
  a_{\theta} \end{pmatrix} = \begin{pmatrix}
  1 & 1 & 1& \cdots &1 \\
  1 & 2 & 2& \cdots &2\\
  1 & 2 & 3& \cdots &3\\
  \vdots  & \vdots & \vdots & \ddots & \vdots\\
  1 & 2 & 3& \cdots & \theta
 \end{pmatrix}
  \begin{pmatrix}
  b_1 \\
  b_2\\
  b_3\\
  \vdots \\
  b_\theta \end{pmatrix}.
  \]
Since the square matrix presented above is clearly invertible, we see that $a_i=b_i$ for every $1\leq i\leq \theta$. From which,  it follows that the elementary representations of $M[f^\infty]$ and $N[f^\infty]$ are the same.  This concludes the proof of the lemma.
\epf

We are now ready to give the proof for Theorems \ref{main thm0} and \ref{main thm0 Selcomplex}.

\bpf[Proof of Theorem \ref{main thm0} and \ref{main thm0 Selcomplex}]
By Lemma \ref{2 Selmer groups} and Remark \ref{strict Sel suffices}, it suffices to show that the assertions of Theorem \ref{main thm0}/\ref{main thm0 Selcomplex} hold for the strict Selmer groups.
For this, we just need to verify that the conditions (I) and (II) specified in Lemma \ref{structure compare} are satisfied for the strict Selmer groups. The validity of (I) follows directly from \cite[Proposition 4.1.3]{LimCMu}. Consequently, all that remains is to show that (II) holds. Let $f$ be an irreducible Weierstrass polynomial, and write $Y=\Op\ps{W}/f^a$ for a given positive integer $a$. Since $\Gal(F_S/F_\cyc)$ acts trivially on $Y$, we have $\Sel_{str}(A/F_\cyc)\ot_{\Op}Y = \Sel_{str}(A\ot_\Op Y/F_\cyc)$ and $\Sel_{str}(A^*/F_\cyc)^\iota\ot_{\Op}Y = \Sel_{str}(A^*\ot_\Op Y^\iota/F_\cyc)^\iota$. Therefore, we are reduced to verifying the following equality
\begin{equation}\label{Gr corank}
   \corank_{\Op} \big(\Sel_{str}(A\ot_\Op Y/F_\cyc)^\Ga\big) = \corank_{\Op} \big(\Sel_{str}(A^*\ot_\Op Y^\iota/F_\cyc)^\Ga\big).
\end{equation}
At this point, it can be readily verified that the conditions \textbf{(C1)-(C4)}, \textbf{(R1)} and \textbf{(R2)} are satisfied for the datum $\big(A\ot_\Op Y, \{A_v\ot_\Op Y\}_{v|p}\big)$. To simplify notation, we write $B=A\ot_\Op Y$, $B_v^-=A_v^-\ot_\Op Y$, $B^*=A^*\ot_\Op Y^\iota$ and $(B_v^*)^-=(A_v^*)^-\ot_\Op Y$. Consider the following commutative diagram
\[ \entrymodifiers={!! <0pt, .8ex>+} \SelectTips{eu}{}\xymatrix{
     0 \ar[r] &  \Sel_{str}(B/F)\ar[r] \ar[d]_{} &  H^1(G_S(F),B) \ar[r] \ar[d]_\al &  \bigoplus_{v\in S} J_v^{str}(B/F) \ar[d]_\be \\
     0 \ar[r] &  \Sel_{str}(B/F_\cyc)^\Ga\ar[r]_{} &   H^1(G_S(F_\cyc),B)^\Ga \ar[r] &  \bigoplus_{v\in S}J_v^{str}(B/F_\cyc)^\Ga
     }\]
with exact rows. From the Hochschild-Serre
spectral sequence, we see that the map $\al$ is surjective with $\ker \al = H^1(\Ga, H^0(F_\cyc, B))$ and
\[ \ker\be = \bigoplus_{v|p} H^1(\Ga_v, H^0(F_{\cyc,v}, B_v^-)). \]
It follows from \textbf{(R1)} and \textbf{(R2)} that these kernels are finite. Therefore, we may apply the snake lemma to conclude that the restriction map
\[ \Sel_{str}(B/F)\lra \Sel_{str}(B/F_\cyc)^\Ga \]
has finite kernel and cokernel. A similar assertion holds for $B^*$. Hence the verification of (\ref{Gr corank}) is reduced to showing the equality
\begin{equation}\label{Gr corank2} \corank_{\Op} \big(\Sel_{str}(B/F)\big) = \corank_{\Op} \big(\Sel_{str}(B^*/F)\big).
\end{equation}
Now, recall that for a cofinitely generated $\Op$-module $N$, we have $|N[\pi^n]| = q^{\corank_{\Op}(N)+O(1)}$, where $q=|\Op/\pi|$. In view of this, the verification of (\ref{Gr corank2}) is then reduced to showing that
\begin{equation}\label{Gr corank3}
 \frac{|\Sel_{str}(B/F)[\pi^n]|}{|\Sel_{str}(B^*/F)[\pi^n]|} = q^{O(1)}.
 \end{equation}
Now, for every $n$, we denote by $\Sel_{str}(B[\pi^n]/F)$ the strict Selmer group of $B[\pi^n]$, where we replace the occurrences of $B$ by $B[\pi^n]$ in the definition of the strict Selmer group. For $C\in \{B, B^*\}$, another control theorem argument shows that the map
\[ \Sel_{str}(C[\pi^n]/F)\lra \Sel_{str}(C/F)[\pi^n] \]
has finite kernel and cokernel which are bounded independent of $n$. Consequently, for the verification of (\ref{Gr corank3}), we just need to show that
\[ \frac{|\Sel_{str}(B[\pi^n]/F)|}{|\Sel_{str}(B^*[\pi^n]/F)|} = q^{O(1)}. \]
By appealing to the argument in \cite[Formula (53)]{G89}, we obtain
\[ \frac{|\Sel_{str}(B[\pi^n]/F)|}{|H^0(G_S(F),B[\pi^n])|}\times
 \prod_{v|p}|H^0(F_{v},B_v^-[\pi^n])| \hspace{2in}\]
 \[ \hspace{2in} = \frac{|\Sel_{str}(B^*[\pi^n]/F)|}
 { |H^0(G_S(F),B^*[\pi^n])|} \times
 \prod_{v|p}|H^0(F_{v}, (B^*_v)^-[\pi^n])| \]
 for every $n$, where we note that the assumption $\textbf{(C4)}$ is crucial in obtaining this said equality.
Since our data $(B,\{B_v\}_{v|p})$ and $(B^*,\{B^*_v\}_{v|p})$ satisfy  $\textbf{(R1)-(R2)}$, the global $H^0$-terms and local $H^0$-terms appearing in the above equality have order $q^{O(1)}$. This in turn proves the validity of (II). With this, the proof of the theorem is therefore completed. \epf

\section{Proof of Theorems \ref{main thm} and \ref{main thm Selcomplex}} \label{Selmer proof}

As a start, we shall assume that $F_\infty$ is a $\Zp^r$-extension of $F$ with $r\geq 2$. We are going to show that the consideration over $F_\infty$ can be reduced to that of $F_\cyc$.
Fix a $\Zp^{r-1}$-extension $L_\infty$ of $F$ contained in $F_\infty$ such that $L_\infty\cap F^\cyc = F$. Writing $H:=\Gal(L_\infty/F)$, there is a natural identification $G\stackrel{\cong}{\lra} \Ga\times H$ which in turn induces isomorphisms
\[ H\cong \Gal(F_\infty/F_\cyc), \quad  \Ga\cong \Gal(F_\infty/L_\infty) \]
of groups and an isomorphism
\[ R\ps{G} \cong (R\ps{H})\ps{\Ga}\]
of commutative profinite rings. Note that $R\ps{H}$ is now isomorphic to a power series ring in $(m+r-1)$-variables.

\bl \label{Selcomplex cyc}
For each $Z\in \{T, T^*\}$, there is a natural isomorphism
\[H^i(SC(Z/F_\infty))\cong H^i(SC(R\ps{H}^\iota\ot_R Z/F_\cyc)).\]
\el

\bpf Since $R\ps{G}^\iota\ot_R T\cong (R\ps{H})\ps{\Ga}^\iota\ot_{R\ps{H}}R\ps{H}^\iota\ot_R T$, we then have
\begin{eqnarray*}
  H^i(SC(T/F_\infty)) &\cong& H^i(SC(R\ps{G}^\iota\ot_R T/F))\quad \mbox{(by \cite[Proposition 8.6.6]{Nek})} \\
   &\cong&  H^i(SC((R\ps{H})\ps{\Ga}^\iota\ot_{R\ps{H}}R\ps{H}^\iota\ot_R T/F))  \\
   &\cong&  H^i(SC(R\ps{H}^\iota\ot_R T/F_\cyc)) \quad \mbox{(by \cite[Proposition 8.6.6]{Nek})}.
\end{eqnarray*}
\epf

Following \cite[8.3.3]{Nek}, we set $F_{R,G}(C) := \Hom_{R,\cts}(R\ps{G}, C)$ for $C\in \{A, A^*\}$. Then it follows from \cite[(8.4.6.6)]{Nek} that there is following ``duality'' diagram
\[  \entrymodifiers={!! <0pt, 0.8ex>+} \SelectTips{eu}{}\xymatrix@R=.6in @C=1in{
     R\ps{H}^{\iota}\ot_{R} T \ar[d]_{-\ot_{R\ps{H}}R\ps{H}^\vee} \ar[r]^{} \ar[dr] &
     R\ps{H} \ot_{R} T^* \ar[l] \ar[dl] \ar[d]^{-\ot_{R\ps{H}}R\ps{H}^\vee}
     \\
    F_{R,H}(A)  \ar[ur]_{}
    &  F_{R,H}(A^*)^\iota \ar[ul]}  \]
of $R\ps{H}$-modules, where the horizontal maps are given by $\Hom_{R\ps{H}}(-, R\ps{H}(1))$ and the two diagonal maps are given by $\Hom_{\cts}(-,\mu_{p^{\infty}})$. It follows from a straightforward verification that the datum
\[\big(R\ps{H}^{\iota}\ot_{R} T, \{R\ps{H}^{\iota}\ot_{R} T_v\}_{v|p}\big)_{R\ps{H},F}\]
satisfies \textbf{(C1)-(C4)}. We now verify that they also satisfy \textbf{(R1)-(R2)} over $R\ps{H}$.

\bl \label{Shapiro R1-R2}
The modules $H^0(F_\cyc, F_{R,H}(A))^\vee$ and $H^0(F_\cyc, F_{R,H}(A^*))^\vee$ are torsion over $R\ps{H}$. For each $v|p$, the modules $H^0(F_{\cyc,v}, F_{R,H}((A^*_v)^-))^\vee$ and $H^0(F_{\cyc,v}, F_{R,H}(A_v^-))^\vee$ are torsion over $R\ps{H}$.
\el

\bpf
Observe that
\begin{equation}\label{H-G}
 \Hom_{R,\cts}\big(R\ps{G}\ot_{R\ps{H}}R\ps{H}, A\big) \cong \Hom_{R\ps{H},\cts}\big(R\ps{G}, \Hom_{R,\cts}(R\ps{H}, A)\big)=F_{R\ps{H},\Ga}\big(F_{R,H}(A)\big).
\end{equation}
From which, we have
\begin{eqnarray*}
  H^0\big(G_S(F_\cyc), F_{R,H}(A)\big) &\cong& H^0\big(G_S(F), F_{R\ps{H},\Ga}\big(F_{R,H}(A)\big)\big)\quad \mbox{(by Shapiro lemma)} \\
   &\cong&  H^0\big(G_S(F), F_{R,G}(A)\big) \quad \mbox{(by (\ref{H-G}))} \\
   &\cong&  H^0\big(G_S(F_\infty), A\big) \quad \mbox{(by Shapiro lemma)}.
\end{eqnarray*}
Since $H\cong \Zp^{r-1}$ with $r-1\geq 1$, the latter is plainly cotorsion over $R\ps{H}$. A similar argument applies for the local cohomology groups.
\epf

We make the following observation.

\br \label{reduction remark}
In the event that  $F_\infty$ is a $\Zp^r$-extension of $F$ such that $F_\cyc\subseteq F_\infty$ and $r\geq 2$,
  Lemmas \ref{Selcomplex cyc} and \ref{Shapiro R1-R2} show that our task for proving the algebraic functional equation for the datum $\big(T, \{T_v\}_{v|p}\big)_{R,F}$ over $F_\infty$  is reduced to proving the algebraic functional equation for the datum $\big(R\ps{H}^{\iota}\ot_{R} T, \{R\ps{H}^{\iota}\ot_{R} T_v\}_{v|p}\big)_{R\ps{H},F}$ over the cyclotomic $\Zp$-extension $F_\cyc$, where the datum $\big(R\ps{H}^{\iota}\ot_{R} T, \{R\ps{H}^{\iota}\ot_{R} T_v\}_{v|p}\big)_{R\ps{H},F}$ is now defined over $R\ps{H}$, a power series ring in $(m+r-1)$-variables.
\er

Taking this observation into account, we will carry out our proof of Theorem \ref{main thm Selcomplex} over $F_\cyc$. The remainder of the section will be devoted to this proof which will proceed by induction on $m$, the number of variables of the power series ring $R=\Op\ps{W_1,...,W_m}$.

We begin with some preliminary  discussion.
Let $\ell$ be a linear element of $R$. Then one has an identification
\[\overline{R}:= R/\ell \cong \Op\ps{W_1,...,W_{m-1}}\]
by Lemma \ref{linear element lemma}.

For a given datum $(T, \{T_v\}_{v|p})_{R,F}$, we obtain another datum $(\overline{T}, \{\overline{T}_v\}_{v|p})_{\overline{R},F}$ which is given by $\overline{T} = T/\ell$ and $\overline{T}_v=T_v/\ell$. One then writes $(\overline{A}, \{\overline{A}_v\}_{v|p})_{\overline{R},F}$ for the discrete datum of  $(\overline{T}, \{\overline{T}_v\}_{v|p})_{\overline{R},F}$. Note that we have $\overline{A} = A[f]$ and $\overline{A}_v=A_v[f]$. It is straightforward to check that if $(T, \{T_v\}_{v|p})_{R,F}$ satisfies  \textbf{(C1)-(C4)}, then $(\overline{T}, \{\overline{T}_v\}_{v|p})_{\overline{R},F}$ will also satisfy these conditions.
However, in order for $(\overline{T}, \{\overline{T}_v\}_{v|p})_{\overline{R},F}$  to satisfy \textbf{(R1)-(R2)}, the linear element $\ell$ has to be chosen appropriately. This is the content of the next lemma.

\bl \label{specialization lemma1}
Suppose that $\ell$ is a linear element of $R$ which does not lie in the support of $A(F_\infty)^\vee$, $A^*(F_\infty)^\vee$, $H^0(F_{\infty,w}, A_v^-)^\vee$ and $H^0(F_{\infty,w}, (A_v^*)^-)^\vee$ for every prime $w$ of $F_\infty$ above $v$. Then the discrete datum $(\overline{A}, \{\overline{A}_v\}_{v|p})$ satisfy \textbf{(R1)-(R2)}.
\el

\bpf
This is immediate from Lemma \ref{coprime torsion}.
\epf

There are natural short exact sequences
\[ 0\lra T \stackrel{\cdot\ell}{\lra} T \lra \overline{T}\lra 0; \]
\[ 0\lra T_v \stackrel{\cdot\ell}{\lra} T_v \lra \overline{T}_v\lra 0. \]
From which, we obtain the following short exact sequence of complexes
\[ 0\lra SC(T/F_\cyc) \stackrel{\cdot\ell}{\lra} SC(T/F_\cyc) \lra SC(\overline{T}/F_\cyc)\lra 0  \]
which in turn yields the following short exact sequence
\[0 \lra H^i(SC(T/F_\cyc))/\ell \lra  H^i(SC(\overline{T}/F_\cyc)) \lra H^{i+1}(SC(T/F_\cyc))[\ell] \lra 0 \] for every $i$.

We can now record the following.

\bl \label{specialization lemma2}
Suppose that $\ell$ is a linear element of $R$ which does not lie in the support of $A(F_\infty)^\vee$, $A^*(F_\infty)^\vee$, $H^0(F_{\infty,w}, A_v^-)^\vee$ and $H^0(F_{\infty,w}, (A_v^*)^-)^\vee$ for every prime $w$ of $F_\infty$ above $v$. Then there is an injection
\[H^2(SC(T/F_\cyc))/\ell \hookrightarrow  H^2(SC(T/F_\cyc)),\]
whose cokernel is a torsion $\overline{R}$-module (and so a pseudo-null $\overline{R}\ps{\Ga}$-module).   There is also an analogous assertion upon replacing all occurrence of $T$ by $T^*$.
\el

\bpf
From the above discussion, we have a short exact sequence
\[0 \lra H^2(SC(T/F_\cyc))/\ell \lra  H^2(SC(\overline{T}/F_\cyc)) \lra H^{3}(SC(T/F_\cyc))[\ell] \lra 0 \]
and an isomorphism
\[ H^{3}(SC(T/F_\cyc))/\ell \cong H^{3}(SC(\overline{T}/F_\cyc)),\]
where we make use of the fact that $H^4(SC(T/F_\cyc))=0$ (see \cite[Proposition 9.7.2(ii) and (iii)]{Nek}).
By Lemma \ref{strict selmer Selcomplex}, we have that $H^3(SC(\overline{T}/F_\cyc))$ is torsion over $\overline{R}$ and $H^3(SC(T/F_\cyc))$ is torsion over $R$. Combining these observations with the above isomorphism, and applying Lemma \ref{coprime torsion}, we deduce that $H^3(SC(T/F_\cyc))[\ell]$ is torsion over $\overline{R}$. The conclusion of the lemma then follows.
\epf

We are ready to give a portion of the proof of our main theorem.

\bpf[Proof of Theorem \ref{main thm Selcomplex}(a) and (b)]
When $m=0$, this is precisely a special case of Theorem \ref{main thm0 Selcomplex}. We may therefore assume that $m\geq 1$ and that Theorem \ref{main thm Selcomplex} holds for $m-1$. We begin with proving (a). Suppose that $H^2(SC(T/F_\cyc))$ is torsion over $R\ps{\Ga}$. Choose a linear element $\ell\in R$ with the properties that $\ell$ does not lie in the support of $A(F_\cyc)^\vee$, $A^*(F_\cyc)^\vee$, $\big(A(F_{\cyc,w})^\vee\big)_{R\mathrm{-tor}}$ and $\big(A^*(F_{\cyc,w})^\vee\big)_{R\mathrm{-tor}}$ for every $w\in S(F_\cyc)$ and that $\ell$ does not lie in the support of $H^2(SC(T/F_\cyc))$. Consequently, we may apply Lemma \ref{coprime torsion} to conclude that $H^2(SC(T/F_\cyc))/\ell$ is torsion over $\overline{R}\ps{\Ga}$, where we write $\overline{R}=R/\ell$. By Lemma \ref{specialization lemma2}, this in turn implies that $H^2(SC(\overline{T}/F_\cyc))$ is torsion over $\overline{R}\ps{\Ga}$. Now since $\overline{R}$ is isomorphic to a power series ring in $(m-1)$-variables by Lemma \ref{linear element lemma}, we may apply our induction hypothesis to conclude that $H^2(SC(\overline{T}^*/F_\cyc))$ is also torsion over $\overline{R}\ps{\Ga}$. From which, it follows from a combination of Lemmas \ref{specialization lemma2} and \ref{coprime torsion} that $X_{str}(A^*/F_\cyc)$ is torsion over $R\ps{\Ga}$, which proves assertion (a) of Theorem \ref{main thm}.

Assertion (b) of Theorem \ref{main thm Selcomplex} has a similar proof as above. \epf

For the final assertion of the main theorem, we require some further preliminary discussion. The following proposition is the first step in this process.

\bp \label{torsion surjective}
Retain the settings as in Theorem \ref{main thm Selcomplex}.
Suppose that $H^2(SC(T/F_\infty))$ is torsion over $R\ps{G}$. Then we have $H^1(SC(T/F_\infty)=0$.   There is also an analogous assertion upon replacing all occurrence of $T$ by $T^*$.
\ep

\bpf
We shall only show that $H^1(SC(T^*/F_\infty))=0$, the case $H^1(SC(T/F_\infty))$ follows by duality.
In view that we have established Theorem \ref{main thm Selcomplex}(a), it follows from the assumption of the proposition that the modules $\Sel_{str}(A/F_\infty)^\vee$, $\Sel_{str}(A^*/F_\infty)^\vee$, $H^2(SC(T/F_\infty))$ and $H^2(SC(T^*/F_\infty))$ are all torsion over $R\ps{G}$.
 By extracting a portion of the long exact sequence (\ref{long selmer seq2}), we have
\begin{multline*}
H^2(G_S(F),A)^\vee \lra H^1(SC(T^*/F))
  \lra\bigoplus_{v|p}H^1(F_v, A_v^-)^\vee\oplus \bigoplus_{v\nmid p} H^1(\Gal(F_v^{ur}/F_v), A^{I_v})^\vee \\ \stackrel{\rho}{\lra} H^1(G_S(F), A)^\vee   \lra H^2(SC(T^*/F)).
\end{multline*}
Since we also have an exact sequence similar to the above for every finite extension $L$ of $F$ contained in $F_\infty$, we may take limit over these extensions and obtain
\begin{multline} \label{long seq3}
H^2(G_S(F_\infty),A)^\vee \lra H^1(SC(T^*/F_\infty))
  \lra\bigoplus_{v\in S}J_v(A/F_\infty)^\vee\\ \stackrel{\rho}{\lra} H^1(G_S(F_\infty), A)^\vee   \lra H^2(SC(T^*/F_\infty))
\end{multline}
Now, recall that the fine Selmer group $R(A^*/F_\infty)$ (in the sense of \cite{CS}; see also \cite{LimFine}) is defined by the following exact sequence
\[0\lra R(A^*/F_\infty) \lra H^1(G_S(F_\infty),A^*) \lra \bigoplus_{v\in S}K_v(A^*/F_\infty),\]
where $K_v(A^*/F_\infty) = \ilim_L \oplus_{w|v}H^1(L_w,A^*)$, and one checks easily that there is a natural injection $R(A^*/F_\infty)\hookrightarrow\Sel_{str}(A^*/F_\infty)$. As a consequence, we see that $R(A^*/F_\infty)^\vee$ is also torsion over $R\ps{G}$. By \cite[Lemma 7.1]{LimFine}, this in turn implies that
\begin{equation} \label{WL}
H^2(G_S(F_\infty),A)=0.
\end{equation}
On the other hand, standard coranks calculations \cite[Propositions 1-3]{G89} tell us that
\[
     \rank_{R\ps{G}} \left(J_v(A/F_\infty)^\vee \right)
    =\left\{
       \begin{array}{ll}
         [F_v:\Qp](d-d_v) , & \hbox{if $v | p$;} \\
         0, & \hbox{otherwise.}
       \end{array}
     \right.
    \]
    and
\[ \rank_{R\ps{G}} \left(H^1(G_S(F_\infty), A)^\vee \right)-\rank_{R\ps{G}} \left(H^2(G_S(F_\infty), A)^\vee\right)
    = r_2(F)d +
 \sum_{v~\mathrm{real}}\left(d-d_v^+)\right),
  \]
  where $r_2(F)$ denotes the number of complex primes of $F$. Putting these observations into the exact sequence (\ref{long seq3}) and taking our standing assumption \textbf{(C4)} into account, we see that $H^1(SC(T^*/F_\infty))$ is torsion over $R\ps{G}$.

To continue, we need to recall some more terminologies from \cite{Nek}. For each $v\in S$, we set
\[ U_v^-(T^*) = \mathrm{Cone}\left(U_v^+(T^*)\stackrel{-i_v}\lra C(\Gal(\overline{F}_v/F_v),T^*)
                \right)
  \]
  and $U_S^-(T^*) = \bigoplus_{v\in S}U_v^-(T^*)$. From \cite[(6.1.3)]{Nek}, we have the following short exact sequence
\[ 0 \lra U^-_S(T^*)  \lra SC(T^*/F_\infty) \lra C(F_\infty, T^*) \lra 0 \]
of complexes, whose cohomology yields an exact sequence
\[ \bigoplus_{v|p}H^0(F_v, (T_v^*)^-) \lra H^1(SC(T^*/F)) \lra H^1(G_S(F), T^*). \]
For every finite extension $L$ of $F$ contained in $F_\infty$, we have a similar exact sequence, and upon taking inverse limit over these extensions, we obtain an injection
\begin{equation}\label{inject1}
   H^1(SC(T^*/F_\infty)) \hookrightarrow H^1_{\mathrm{Iw}}(F_\infty, T^*),
\end{equation}
where we write $H^1_{\mathrm{Iw}}(F_\infty, T^*):=\plim_L H^1(G_S(L), T^*)$. Now consider the following spectral sequence (see \cite[Theorem 8.5.6]{Nek})
\[ \Ext_{R\ps{G}}^i\left((H^j(G_S(F_\infty), A^*)^\vee, R\ps{G}\right) \Rightarrow H^{i+j}_{\mathrm{Iw}}(F_\infty, T^*), \]
whose low-degree terms gives the following exact sequence
\[ 0\lra\Ext_{R\ps{G}}^1\left(H^0(G_S(F_\infty), A^*)^\vee, R\ps{G}\right) \lra H^1_{\mathrm{Iw}}(F_\infty, T^*)\lra \Ext_{R\ps{G}}^0\left(H^1(G_S(F_\infty), A^*)^\vee, R\ps{G}\right). \]
 Since the module $H^0\left(G_S(F_\infty),  A^*)\right)^\vee$ is torsion over $R$ by our hypothesis \textbf{(R1)}, it is pseudo-null over $R\ps{G}$ for the dimension of the group $G$ is $\geq 1$. It then follows that the $\Ext^1$-term in the above exact sequence vanishes. From which, we obtain an inclusion
 \[ H^1_{\mathrm{Iw}}(F_\infty, T^*)\hookrightarrow \Ext_{R\ps{G}}^0\left(H^1(G_S(F_\infty), A^*)^\vee, R\ps{G}\right)\] which upon combined with (\ref{inject1}) yields an inclusion
 \[ H^1(SC(T^*/F_\infty))\hookrightarrow \Ext_{R\ps{G}}^0\left(H^1(G_S(F_\infty), A^*)^\vee, R\ps{G}\right). \]
 But an $\Ext_{R\ps{G}}^0$-term does not contain non-trivial torsion $R\ps{G}$-submodule, and so this forces $H^1(SC(T^*/F_\infty))$ to vanish in view of the torsionness observation in the previous paragraph.
\epf

We require one more lemma.

\bl \label{specialization lemma3}
Suppose that $H^2(SC(T/F_\cyc))$ is torsion over $R\ps{\Ga}$.
 For every linear element $\ell$ of $R$ which does not lie in the support of $H^2(SC(T/F_\cyc))$, $A(F_\cyc)^\vee$ and $\big(A_w^-(F_{\cyc,w})^\vee\big)_{R\mathrm{-tor}}$ for every $w\in S(F_\cyc)$, we have $\ell \in \mathcal{L}_{R\ps{\Ga}}\big(H^2(SC(T/F_\cyc))\big)$.
\el

\bpf
For the verification of this lemma, we shall appeal to Lemma \ref{linear element choice}.
Our choice of $\ell$ guarantees that $H^2(SC(T/F_\cyc))/\ell$ is torsion over $\overline{R}\ps{\Ga}$. Therefore, it suffices to show that $H^2(SC(T/F_\cyc)_{\mathrm{null}}[\ell]$ is pseudo-null over $\overline{R}\ps{\Ga}$. We shall, in fact, show that $H^2(SC(T/F_\cyc)[\ell] =0$ which will imply what we intend to establish.
By Lemma \ref{specialization lemma2}, we see that $H^2(SC(\overline{T}/F_\cyc))$ is also torsion over $\overline{R}\ps{\Ga}$, and this in turn implies that $H^1(SC(\overline{T}/F_\cyc))=0$ by the preceding proposition. The required vanishing then follows from this and the following short exact sequence
\[0 \lra H^1(SC(T/F_\cyc))/\ell \lra  H^1(SC(\overline{T}/F_\cyc)) \lra H^2(SC(T/F_\cyc))[\ell] \lra 0. \]
\epf

We can now complete the proof of Theorem \ref{main thm Selcomplex}.

\bpf[Proof of Theorem \ref{main thm Selcomplex}(c)]
Since the support of $A(F_\cyc)^\vee$, $A^*(F_\cyc)^\vee$, $A(F_{\cyc,w})^\vee$, $A^*(F_{\cyc,w})^\vee$, $H^2(SC(T/F_\cyc))$ and $H^2(SC(T^*/F_\cyc))$ only contains finitely many primes of height one, we can always find an infinite collection of linear ideals (which we denote by $\{\ell_i\}$) that do not lie in the support of these modules. By Lemma \ref{specialization lemma3}, these linear ideals lie in \[\mathcal{L}_{R\ps{G}}\big(H^2(SC(T/F_\cyc))\big)\cap \mathcal{L}_{R\ps{G}}\big(H^2(SC(T^*/F_\cyc))\big).\]
For each $i$, we denote by $\{T_i, T_{i,v})_{v|p}\}_{R_i,F}$ the datum given by $T_i:=T/\ell_i$ and $T_{i,v}:=T_v/\ell_i$, where $R_i:=R/\ell_i$. By our choice of linear ideals $(\ell_i)$'s, we may apply Lemma \ref{specialization lemma1} to conclude that the datum $\{T_i, T_{i,v})_{v|p}\}_{R_i,F}$ satisfies \textbf{(C1)-(C4)}, \textbf{(R1)} and \textbf{(R2)}. On the other hand, it follows from Lemma \ref{specialization lemma2} that $H^2(SC(Z/F_\cyc))/\ell_i$ is torsion over $R_i\ps{\Ga}$ with
\[ \ch_{R_i\ps{\Ga}} \Big(\big(H^2(SC(Z/F_\cyc))/\ell_i\Big) =\ch_{R_i\ps{\Ga}} \Big(H^2(SC(Z_i/F_\cyc))\Big) \]
for $Z\in \{T,T^*\}$. On the other hand, in view of Lemma \ref{linear element lemma}, we may apply our induction hypothesis to obtain
\[ \ch_{R_i\ps{\Ga}} \Big(H^2(SC(T_i/F_\cyc)) \Big) =\ch_{R_i\ps{\Ga}} \Big(H^2(SC(T_i^*/F_\cyc))\Big). \]
We may therefore apply Proposition \ref{alg lemma} to obtain the conclusion of the theorem.
\epf

\section{Examples} \label{examples section}

\subsection{Abelian variety with good ordinary reduction at every $p$-adic prime}

Let $\mathcal{A}$ be a $g$-dimensional abelian variety defined over a number field $F$, with the property that $\mathcal{A}$ has good ordinary reduction at every prime of $F$ above $p$. In this context, we set $T$ as the Tate module of $\mathcal{A}$ which is a free $\Zp$-module of rank $2g$. By appealing to \cite[Proposition 1.1]{GH},
we see that $T_v^+$ is a free $\Zp$-module of rank $d$ for every real prime $v$. For each $v|p$, it follows from the discussion in \cite[P150-151]{CG} that there is $\Gal(\bar{F_v}/F_v)$-submodule, which will take the role of our $A_v$, of $\mathcal{A}[p^\infty]$ such that the inertia group acts trivially on $A_v^-:= \mathcal{A}[p^\infty]/A_v$. Moreover, there is a natural identification $A_v^- \cong \widetilde{\mathcal{A}}[p^\infty]$, where $\widetilde{\mathcal{A}}$ is the mod-$p$ reduction of $\mathcal{A}$, which is an abelian variety over the residue field of $F_v$. We set $T_v$ to be the Tate module of $A_v$. Taking into account all of the above observations, one can check that the datum $(T, (T_{v})_{v|p})_{\Z, F}$ satisfies \textbf{(C1)-(C4)}. On the other hand, property \textbf{(R1)} follows as a consequence of Imai's theorem \cite{Imai}. For \textbf{(R2)}, we first observe that $H^0(F_{\cyc,w}, A_v^-)$ identifies with the subgroup of $\widetilde{\mathcal{A}}[p^\infty]$ consisting of points in the residue field of $F_{\cyc,w}$. The latter is plainly a finite field in view that the inertia group of $F_{\cyc,w}/F_v$ is a subgroup of $\Gal(F_{\cyc,w}/F_v)$ with finite index. Hence it follows that $H^0(F_{\cyc,w}, A_v^-)$  is finite. Therefore, Theorem \ref{main thm0}/\ref{main thm0 Selcomplex} and Theorem \ref{main thm}/\ref{main thm Selcomplex} apply. We would like to point out that Lee \cite{Lee} has already proven this result in the situation of a cyclotomic $\Zp$-extension. %However, we also like to highlight that our theorem are applicable to multiple $\Zp$-extensions$-$a context that Lee's result does not cover.

We now say something about the torsionness of Selmer groups. For a $\Zp^r$-extension $F_\infty$ of $F$ containing $F_\cyc$, a classical conjecture of Mazur \cite{Maz} (also see \cite[Conjecture 7.4]{OV})  posits that $X_{Gr}(\mathcal{A}/F_\infty)$ is a torsion $\Zp\ps{G}$-module, where $G=\Gal(F_\infty/F)$. When $\mathcal{A}$ is an elliptic curve over $\Q$ and $F$ is an abelian extension of $\Q$, a deep result of Kato \cite{K} says that $X_{Gr}(\mathcal{A}/F_\cyc)$ is a torsion $\Zp\ps{\Gal(F_\cyc/F)}$-module. Combining this with a result of Hachimori-Ochiai \cite{HO}, one can deduce that $X_{Gr}(\mathcal{A}/F_\infty)$ is a torsion $\Zp\ps{\Gal(F_\cyc/F)}$-module for every $\Zp^r$-extension $F_\infty$ of $F$ containing $F_\cyc$.

We should mention that the torsionness conjecture is known to hold, thanks to a deep result of Kato \cite{K}, in the case where $A$ is an elliptic curve over $\Q$ and $F$ is an abelian extension of $\Q$.

\subsection{$p$-ordinary modular forms}

Let $f=\sum_{n\geq 1}a_nq^n$ be a normalized new cuspidal modular eigenform of weight $k\geq 2$, level $N$ and nebentypus $\epsilon$. We shall always assume that the coefficient $a_p$ is a $p$-adic unit (in other words, the modular form $f$ is $p$-ordinary). In the event that $k=2$, we assume further that the (odd) prime $p$ does not divide the level $N$. (This is to ensure that the absolute value of the unit root of $f$ is not one which is required for the verification of \textbf{(R2)}; see discussion below.) Let $K_f$ be the number field obtained by adjoining all the Fourier coefficients of $f$ and values of $\epsilon$ to $\Q$. Throughout our discussion, we fix a prime $\mathfrak{p}$ of $K_f$ above $p$. We then let $V_f$ denote the corresponding two-dimensional $K_{f,\mathfrak{p}}$-linear Galois representation attached to $f$ in the sense of Deligne. Writing $\Op=\Op_{K_{f,\mathfrak{p}}}$ for the ring of integers of $K_{f,\mathfrak{p}}$, we denote by $T_f$ the $\Gal(\bar{\Q}/\Q)$-stable $\Op$-lattice in $V_f$ as defined in \cite[Section 8.3]{K}. By our $p$-ordinarity assumption, there is a short exact sequence
\[ 0 \lra T_{f,p} \lra T_f \lra T_{f,p}^-\lra 0\]
of $\Gal(\overline{\Q}_p/\Qp)$-modules, where $T_{f,p}$ and $T_{f,p}^-$ are free $\Op$-modules of rank $1$. Furthermore, the Galois group $\Gal(\overline{\Q}_p/\Qp)$ acts on $T_{f,p}$ via an unramified character $\delta_f$ with the property that $\delta_f$ sends the (geometric) Frobenius at $p$ to the $p$-adic unit to the unit
root of the Hecke polynomial of $f$ at $p$. Additionally, the same Galois group $\Gal(\overline{\Q}_p/\Qp)$ acts on $T_{f,p}^-$ via $\delta_f^{-1}\epsilon\kappa^{1-k}$, where $\kappa$ denotes the $p$-adic cyclotomic character.

 We shall proceed to verify that the datum $(T_f, T_{f,p})$ satisfies \textbf{(R1)} and \textbf{(R2)}, leaving the relatively straightforward verification of \textbf{(C1)-(C4)} to the readers. For this, we work with the discrete datum $(A, A_p)$, where $A_f:= V_f/T_f = T_f \ot K_f/\Op$ and $A_{f,p}: = T_{f,p} \ot K_f/\Op$.

We now present the following general lemma which will be required for the subsequent discussion.

\bl \label{irreducible R1}
Let $F$ be either a number field or a finite extension of $\Q_p$. Let  $\rho :\Gal(\overline{F}/F)\lra \mathrm{GL}(V)$ be a Galois representation, where $V$ is finite dimensional $K$-vector space for some finite extension $K$ of $\Qp$. Suppose that $\rho$ is irreducible and that $\dim (\mathrm{Im}~\rho)>1$, where here $\dim$ means the dimension as a $p$-adic Lie group. Then for every finite extension $L$ of $F$, we have $H^0(\Gal(\overline{F}/L_\cyc), V)=0$.
\el

\bpf
Clearly, we may assume that $L$ is Galois over $F$, and from this assumption, it follows that $\Gal(\overline{F}/L_\cyc)$ is a normal subgroup of $\Gal(\overline{F}/F)$. Suppose on the contrary that $V_0:=H^0(\Gal(\overline{F}/L_\cyc), V)\neq 0$. By normality, we see that $V_0$ is a $\Gal(\overline{F}/F)$-subrepresentation of $V$. In view of the irreducibility of $V$, we must therefore have $V_0 =V$. Consequently, the map $\rho$ factor through $\Gal(L_\cyc/F)$ to yield a map
\[ \overline{\rho}: \Gal(L_\cyc/F) \lra \mathrm{GL}(V), \]
whose image coincides with $\mathrm{Im}~\rho$. However, since the group $\Gal(L_\cyc/F)$ can clearly be expressed as a semi-direct product of a copy of $\Zp$ with a finite group, its image under $\overline{\rho}$ cannot be a Lie group of dimension $\geq 2$. Thus, we arrive at a contradiction, which proves the lemma.
\epf

Equipped with the aforementioned lemma, we are now in a position to verify \textbf{(R1)}. Indeed, it is well-known that $V_f$ is irreducible as a $\Gal(\overline{\Q}/\Q)$-module. Furthermore, the dimension of the image of the map
\[\rho_f: \Gal(\overline{\Q}/\Q)\lra \mathrm{GL}(V_f)\]
is $\geq 2$ (see \cite{Rib}). Hence it follows from Lemma \ref{irreducible R1} that $H^0(\Gal(\overline{\Q}/\Q_\cyc), V_f)=0$, which is equivalent to saying that $H^0(\Gal(\overline{\Q}/\Q_\cyc), A_f)$ is finite. For the dual datum, one simply notes that $V_f^{*}=\Hom_{K_{f,\p}}(V_f, K_{f,\p}(1))$ and that taking dual preserves irreducibility and the $p$-adic dimension of the image. Hence this shows \textbf{(R1)}.

We come to verifying \textbf{(R2)}. For this, we need to return to the short exact sequence
\[ 0 \lra T_{f,p} \lra T_f \lra T_{f,p}^-\lra 0\]
of $\Gal(\overline{\Q}_p/\Qp)$-modules, and utilize certain description of the modules. In particular, $T_{f,p}$ is the unramified subrepresentation with character $\delta_f$, and $T_{f,p}^-$ is the ramified quotient with character $\delta_f^{-1}\epsilon\kappa^{1-k}$, where $\kappa$ is the $p$-adic cyclotomic character. For the dual datum $(T_f^*, (T_{f,p}^-)^*)$, we then have
\[ 0\lra  (T_{f,p}^-)^*\lra T_f^* \lra T_{f,p}^* \lra 0,\]
where $T_{f,p}^*$ comes equipped with the character $\delta_f^{-1}\kappa$. Since $T_{f,p}^-$ and $T_{f,p}^*$ are free $\Op$-modules of rank $1$, the validity of \textbf{(R2)} will follow once we can show that the characters $\delta_f^{-1}\epsilon\kappa^{1-k}$ and $\delta_f^{-1}\kappa$, upon restricted to $\Gal(\overline{\Qp}/\Q_{p,\cyc})$, remain non-trivial. We shall show the slightly stronger assertion, namely, $\delta_f^{-1}\epsilon\kappa^{1-k}|_{\Gal(\overline{\Qp}/\Qp(\mu_{p^\infty}))}$ and $\delta_f^{-1}\kappa|_{\Gal(\overline{\Qp}/\Qp(\mu_{p^\infty}))}$ are non-trivial.  Since $\delta_f$ is unramified, upon restriction, we have $\delta_f^{-1}\epsilon\kappa^{1-k}|_{\Gal(\overline{\Qp}/\Qp(\mu_{p^\infty}))} = \delta_f^{-1}\epsilon$ and
$\delta_f^{-1}\kappa|_{\Gal(\overline{\Qp}/\Qp(\mu_{p^\infty}))}= \delta_f^{-1}$. But as seen above, $\delta_f$ sends the Frobenius at $p$ to the unit
root of the Hecke polynomial of $f$ at $p$, and this unit root cannot be equal to $1$ for it has complex absolute value $p^{(k-1)/2}$. Therefore, $\delta_f$ is non-trivial, and in fact, cannot be of finite order. Consequently, $\delta_f^{-1}\epsilon$ is also non-trivial. In conclusion, we have that the datum $(T_f, T_{f,p})$ satisfies \textbf{(C1)-(C4)}, \textbf{(R1)} and \textbf{(R2)}.

Now, let $F$ be a finite Galois extension of $\Q$, and consider the datum $(T_f, (T_{f,p})_{v|p})_{\Op,F}$ obtained from  by base changing to $F$. Since $F$ is Galois, it is either totally real or totally imaginary. Therefore, one may apply Lemma \ref{data base change} to conclude that the base-changed datum still satisfy  \textbf{(C1)-(C4)}. The validity of \textbf{(R1)} still follows from Lemma \ref{irreducible R1}. For \textbf{(R2)}, one simply notes that $\delta_f$ is an unramified character of infinite order, and so for every prime $w$ of $F_\cyc$ above $v$, we see that upon restricting to $F_{\cyc,w}$, $\delta_f$ remains a character with infinite order. Therefore, the same argument as above carries over.

Hence Theorem \ref{main thm0}/\ref{main thm0 Selcomplex} and Theorem \ref{main thm}/\ref{main thm Selcomplex} are applicable, and we can derive a previous result of Jha-Pal \cite[Theorem 3.8]{JP} from our Theorem. It's worth mentioning that we also obtain a corresponding result for multiple $\Zp$-extensions, a context not addressed in the work of Jha-Pal. We would also like to note that, through an analysis akin to the one presented above, our main result applies to yield algebraic functional equation results for a nearly ordinary Hilbert modular form, which not only recovers but also extends the prior work of Jha-Majumdar \cite{JM} in this particular situation.

Finally, we say a bit on the torsionness of the Selmer groups. By a deep result of Kato \cite{K}, the torsionness of $\Sel(A_f/F_\cyc)^\vee$ is known for an elliptic cuspidal modular eigenform $f$ when $F$ is an abelian extension of $\Q$. Building on this and following a similar argument to that in \cite[Theorem 2.3]{HO}, one can see that $\Sel(A_f/F_\cyc)^\vee$ is torsion over $\Zp\ps{\Gal(F_\infty/F)}$ for every $\Zp^r$-extension $F_\infty$ of $F$ which contains $F_\cyc$. Therefore, in this context, we do obtain algebraic functional equation results unconditionally.

\subsection{Rankin-Selberg products}

Let $f_i$ ($i=1,2$) be two normalized new cuspidal modular eigenforms of weight $k_i \geq 2$, level $N_i$ and nebentypus $\epsilon_i$. In this subsection, we always assume that $k_1>k_2$, and that $f_1$ is $p$-ordinary and does not have complex multiplication. We also assumed that $p$ does not divide $N_1N_2$. Let $\mathcal{K}$ be the number field obtained by adjoining all the Fourier coefficients of $f_1, f_2$ and values of $\epsilon_1, \epsilon_2$ to $\Q$. As before, fix a prime $\mathfrak{p}$ of $\mathcal{K}$ above $p$. Writing $\Op=\Op_{\mathcal{K}_{\mathfrak{p}}}$ for the ring of integers of the local field $\mathcal{K}_{\mathfrak{p}}$, we let $V_i$ denote the corresponding two-dimensional $K_{\mathfrak{p}}$-linear Galois representation attached to $f_i$, and $T_{i}$ the corresponding $\Gal(\bar{\Q}/\Q)$-stable $\Op$-lattice.
 By $p$-ordinarity, we have a short exact sequence
\[ 0 \lra T_{1,p} \lra T_{1} \lra T_{1,p}^-\lra 0\]
of $\Gal(\overline{\Q}_p/\Qp)$-modules, where $T_{1,p}$ and $T_{1,p}^-$ are free $\Op$-modules of rank $1$ with $\Gal(\overline{\Q}_p/\Qp)$ acting on $T_{1,p}$ via an unramified character $\delta_{1}$ of infinite order, and on $T_{1,p}^-$ via $\delta_1^{-1}\epsilon_1\kappa^{1-k_1}$. Here, $\kappa$ denotes the $p$-adic cyclotomic character, as previously defined. In this context, our datum $(T, T_p)_{\Op,\Q}$ is defined by $T:=T_{1}\ot T_{2}$ and $T_p:=T_{1,p}\ot T_{2}$.
Again, it can be verified directly that this datum satisfies \textbf{(C1)-(C4)}.

As will be shown in Section \ref{appendix}, particularly in Proposition \ref{product modular forms}, the representation $V:= V_1\ot_{\mathcal{K}_{\mathfrak{p}}} V_2$ is irreducible. Once we have this irreducibility, we may proceed as in the case of one modular form by appealing to Lemma \ref{irreducible R1} to obtain the validity of \textbf{(R1)}.

We now come to \textbf{(R2)}. For this, we break into two cases: $
f_2$ is $p$-ordinary or not.
Suppose first that $f_2$ is non-ordinary at $p$. Then in this case, $V_2$ is an irreducible $\Gal(\overline{\Qp}/\Qp)$-module. Therefore, if we let $V_{1,p}^-$ denote the one-dimensional vector space associated to $T_{1,p}^-$, the tensor product $V_{1,p}^-\ot V_2$ remains irreducible over $\Gal(\overline{\Qp}/\Qp)$. The validity of \textbf{(R2)} then follows from an application of Lemma \ref{irreducible R1}.

We next verify \textbf{(R2)} for the case where $f_2$ is ordinary at $p$. Under this assumption, we have a short exact sequence
\[ 0 \lra V_{2,p} \lra V_{2} \lra V_{2,p}^-\lra 0\]
of $\Gal(\overline{\Q}_p/\Qp)$-modules with $\Gal(\overline{\Q}_p/\Qp)$ acting on $V_{2,p}$ via an unramified character $\delta_{2}$ of infinite order, and on $V_{2,p}^-$ via $\delta_2^{-1}\epsilon\kappa^{1-k_2}$. From which, we have  $\Gal(\overline{\Q}_p/\Qp)$ acting on $V_{1,p}^-\ot V_{2,p}$ via $\delta_1^{-1}\delta_{2}\epsilon_1\kappa^{1-k_1}$, and on $V_{1,p}^-\ot V_{2,p}^-$ via $\delta_1^{-1}\delta_2^{-1}\epsilon_1\epsilon_2\kappa^{2-k_1-k_2}$. Now, the image of Frobenius at $p$ under $\delta_1^{-1}\delta_{2}$ (resp. $\delta_1^{-1}\delta_{2}^{-1}$) has complex absolute value $p^{(k_1-k_2)/2}$
(resp., $p^{(2-k_1-k_2)/2}$). Therefore, we may follow the same argument as in the preceding section to verify that \textbf{(R2)} holds. Therefore, our main result applies to yield algebraic functional equation results for Rankin-Selberg products.

\subsection{Hida family} \label{Hida family subsection}

Let $\mathcal{F}$ be an ordinary cuspidal Hida family defined over $\Q$ (in the sense of \cite{Hida86}). Denote by $R_\mathcal{F}$ the integral closures of the irreducible components of the ordinary primitive cuspidal Hecke algebras associated to $\mathcal{F}$. In general, this is a finite integral
extension of the power series ring $\Zp\ps{W}$. We shall assume that $R_\mathcal{F}$ is a power series ring for our subsequent discussion (see \cite{Gou92} for sufficient conditions for ensuring this property).
We also assume that the residual representations associated to $\mathcal{F}$ is absolutely irreducible, and that the restriction, to the decomposition subgroup $\Gal(\overline{\Qp}/\Qp)$, of this
residual representation has non-scalar semi-simplifications.

Let $S$ be a finite set of primes of $\Q$ which is assumed to contain $p$, the primes dividing the tame level of $\mathcal{F}$  and $\infty$. Associated to $\mathcal{F}$, we have a $G_S(\Q)$-module $T_\mathcal{F}$ which is a free $R_\mathcal{F}$-module of rank 2 and sits in the following short exact sequence
\begin{equation}\label{local Hida}
  0\lra T_{\mathcal{F},p} \lra T_\mathcal{F}\lra T_{\mathcal{F},p}^-\lra 0
\end{equation}
of $\Gal(\overline{\Qp}/\Qp)$-modules. Our datum in this context is $(T_\mathcal{F}, T_{\mathcal{F},p})_{ R_\mathcal{F},\Q}$, whose discrete datum is in turn denoted by $(A_\mathcal{F}, A_{\mathcal{F},p})_{ R_\mathcal{F},\Q}$. We also write $A^-_{\mathcal{F},p} = A_\mathcal{F}/A_{\mathcal{F},p} = T_{\mathcal{F},p}^-\ot_{R_\mathcal{F}}R_\mathcal{F}^\vee$.

Let $F$ be a finite extension of $\Q$. Suppose that $\p_f$ is a prime ideal of $R_\mathcal{F}$ such that $T_f:=T/\p_f$ corresponds to the Galois representation of a cuspidal modular form $f$ of weight $>2$. This fit into a short exact sequence
\[ 0\lra T_{f,p} \lra T_f \lra T_{f,p}^- \lra 0\]
of $\Op[\Gal(\overline{\Q}_p/\Qp)]$-modules with $T_{f,p}= T_p/\p_f$ and $T_{f,p}^-= T_p^-/\p_f$ being free $\Op$-modules of rank one. If we let $(A_f, A_{f,p})$ be the corresponding discrete datum of $(T_f, T_{f,p})$, one can check easily that $A[\p_f] =A_f$ and $A^-_p[\p_f] =A^-_{f,p}$. Then as seen in previous subsection, we have that $H^0(G_S(F_\cyc), A_f)$ is finite. Since $H^0(G_S(F_\cyc),A_\mathcal{F})[\p_f] = H^0(G_S(F_\cyc), A_f)$, we may apply Lemma \ref{coprime torsion} to deduce that $H^0(G_S(F_\cyc),A_\mathcal{F})$ is cotorsion over $R_\mathcal{F}$. A similar argument shows that $H^0(G_S(F_\cyc),A_\mathcal{F}^*)$ is cotorsion over $R_\mathcal{F}$, and this shows that $\textbf{(R1)}$ is valid. The validity of $\textbf{(R2)}$ can be established through an analogous argument.

We can also consider the product of two Hida families  $\mathcal{F}$ and  $\mathcal{G}$. In this situation, we have two choices of datum which will take the form of $(T_\mathcal{F}\ot T_\mathcal{G}, T_{\mathcal{F},p}\ot T_\mathcal{G})_{ R_\mathcal{F}\widehat{\ot} R_\mathcal{G},\Q}$ or $(T_\mathcal{F}\ot T_\mathcal{G}, T_{\mathcal{F}}\ot T_{\mathcal{G},p})_{ R_\mathcal{F}\widehat{\ot} R_\mathcal{G},\Q}$ accordingly to the choice of the dominant Hida family. where $T_{\mathcal{F},p}$ is defined as in (\ref{local Hida}). One can verify directly both data satisfy \textbf{(C1)-(C4)}. Properties $\textbf{(R1)}$ and $\textbf{(R2)}$ can be shown by a similar argument as above, but this time specialized to the product of two modular forms. Therefore, we obtain an algebraic functional equation for the characteristic ideals of the Selmer groups which aligns with expectations since these ideals are conjecturally related to the $p$-adic $L$-functions constructed by Hida in \cite{Hida85, Hida88} (also see \cite{CH}).

\subsection{Triple product of Hida families}

Let $\mathcal{F}$, $\mathcal{G}$ and $\mathcal{H}$ be three ordinary cuspidal Hida families defined over $\Q$ with coefficient ring $R_\mathcal{F}=R_\mathcal{G}=R_\mathcal{H}=\Op\ps{W}$. Write $T_{\mathcal{F},\mathcal{G},\mathcal{H}}= T_\mathcal{F}\ot T_\mathcal{G}\ot T_\mathcal{H}$ and $R_{\mathcal{F},\mathcal{G},\mathcal{H}}= R_\mathcal{F}\ot R_\mathcal{G}\ot R_\mathcal{H}$. We also assume that the three Hida families $\mathcal{F}$, $\mathcal{G}$ and $\mathcal{H}$ have no complex multiplication. By the ordinarity assumption, there is a short exact sequence
\begin{equation*}
  0\lra T_{\mathcal{J},p} \lra T_\mathcal{J}\lra T_{\mathcal{J},p}^-\lra 0
\end{equation*}
of $\Gal(\overline{\Q}_p/\Qp)$-modules for each $\mathcal{J}\in\{\mathcal{F},\mathcal{G},\mathcal{H}\}$. Then there are three natural choices of datum:
\begin{eqnarray*}
 (T_{\mathcal{F},\mathcal{G},\mathcal{H}}, T_{\mathcal{F},p}\ot T_\mathcal{G}\ot T_\mathcal{H})_{R_{\mathcal{F},\mathcal{G},\mathcal{H}},\Q}, \\
 (T_{\mathcal{F},\mathcal{G},\mathcal{H}}, T_{\mathcal{F}}\ot T_{\mathcal{G},p}\ot T_\mathcal{H})_{R_{\mathcal{F},\mathcal{G},\mathcal{H}},\Q}, \\
  (T_{\mathcal{F},\mathcal{G},\mathcal{H}}, T_{\mathcal{F}}\ot T_\mathcal{G}\ot T_{\mathcal{H},p})_{R_{\mathcal{F},\mathcal{G},\mathcal{H}},\Q},
\end{eqnarray*}
depending on the choice of the dominant Hida family. In these situations, the characteristic ideal of the Selmer groups over $\Q_\cyc$ should be the algebraic counterpart of the ``unbalanced''  $4$-variable $p$-adic $L$-functions which seems to be conjecturally at this point of writing (Nevertheless, one can see the work of Hsieh \cite{Hsi} for a construction of a $3$-variable $p$-adic $L$-function in this direction but without the cyclotomic variable).

Besides the above three data, there is also another datum of interest, whose Selmer group should conjecturally be related to the ``balanced'' $4$-variable $p$-adic $L$-function as constructed recently by Hsieh-Yamana \cite{HsiY}. In this context, the datum is given by $(T_{\mathcal{F},\mathcal{G},\mathcal{H}}, T_{\mathrm{bal},p})_{R_{\mathcal{F},\mathcal{G},\mathcal{H}},\Q}$, where
\[ T_{\mathrm{bal},p} := T_{\mathcal{F},p}\ot T_{\mathcal{G},p}\ot T_\mathcal{H} + T_{\mathcal{F},p}\ot {T_\mathcal{G}}\ot T_{\mathcal{H},p} + T_{\mathcal{F}}\ot T_{\mathcal{G},p}\ot T_{\mathcal{H},p}.\]
Note that $T_{\mathrm{bal},p}$ is a free $R$-module of rank $4$.
From this, it is straightforward to verify that these data satisfy \textbf{(C1)-(C4)}. For $\textbf{(R1)}$, we choose appropriate specializations $f,g,h$ of $\mathcal{F}$, $\mathcal{G}$, $\mathcal{H}$ so that Proposition \ref{product three modular forms} applies to yield the irreduciblity of $V_f\ot V_g\ot V_h$ which shows that $V_f\ot V_g\ot V_h$ satisfies $\textbf{(R1)}$. From which, a similar argument to that in Subsection \ref{Hida family subsection} via Lemma \ref{coprime torsion corollary} will establish $\textbf{(R1)}$. For $\textbf{(R2)}$, it remains to show that $H^0(\Q_{\cyc,p} X)$ is cotorsion over $R_{\mathcal{F},\mathcal{G},\mathcal{H}}$, where
\[ X\in \{T_{\mathcal{F},p}\ot T_{\mathcal{G},p}^-\ot T_{\mathcal{H},p}^-, T_{\mathcal{F},p}^-\ot T_{\mathcal{G},p}\ot T_{\mathcal{H},p}^-, T_{\mathcal{F},p}\ot T_{\mathcal{G},p}^-\ot T_{\mathcal{H},p}^-, T_{\mathcal{F},p}^-\ot T_{\mathcal{G},p}^-\ot T_{\mathcal{H},p}^-\}. \]
For the verification of these, one chooses specializations appropriately such that the character acting on the specialization contains a component in the form of an unramified character which sends the Frobenius at $p$ to a complex value with absolute value not equal to 1. This will yield $\textbf{(R2)}$ for this specialized module. Subsequently, we may apply Lemma \ref{coprime torsion corollary} to obtain $\textbf{(R2)}$ for our data.

\subsection{Half-ordinary Rankin-Selberg universal deformation}

This class of examples is taken from the recent work of the first named author and Loeffler \cite{HL}. Let $\mathcal{F}$ be an ordinary cuspidal Hida families defined over $\Q$ with coefficient ring $R_\mathcal{F}$ which is assumed to be a power series ring. Let $T_{\mathrm{univ}}$ be the 2-dimensional universal deformation lattice with coefficient ring $R_{\mathrm{univ}}$. We shall assumed that  $R_{\mathrm{univ}}\cong \Op\ps{X_1, X_2, X_3}$ (see \cite{Bos,Wes} for discussions on sufficient conditions for ensuring this property). We also assumed that the residual representations of both $T_{\mathcal{F}}$ and  $T_{\mathrm{univ}}$  are absolutely irreducible.
The half-ordinary Rankin-Selberg universal deformation (in the sense of \cite{HL, Loe23}) is the Galois lattice $T_{\mathrm{HL}}:=T_\mathcal{F}\widehat{\ot}_\Op T_{\mathrm{univ}}^*$, where $T_{\mathrm{univ}}^*= \Hom_{R_\mathrm{univ}}(T_{\mathrm{univ}}, R_{\mathrm{univ}}(1))$. In \cite{HL}, Hao-Loeffler constructed a $4$-variables $p$-adic $L$-function associated to this half ordinary universal deformation and showed that the $p$-adic $L$-function satisfied a functional equation.

An initial naive thought is to perhaps consider the Selmer complex of the datum $(T_{\mathrm{HL}}, T_{\mathcal{F},p}\widehat{\ot}_\Op T_{\mathrm{univ}}^*)_{R,\Q}$ over the cyclotomic $\Zp$-extension $\Q_\cyc$. However, such a resulting Selmer complex's cohomology groups will have module structures defined over a power series ring with $5$-variables which definitely does not align with the analytic situation. The correct approach (in the sense of matching the analytic picture) is as follows: write $T_{\mathrm{tame}}$ for the deformation lattice with tame determinant, and $R_{\mathrm{tame}}$ for its coefficient rings (in the sense of \cite[Subsection 4.2]{HL}). It then follows from \cite[Lemma 4.1]{HL}  that $T_{\mathrm{univ}} \cong T_{\mathrm{tame}}\widehat{\ot}_\Op \Op\ps{\Ga}$ with $R_{\mathrm{univ}} \cong R_{\mathrm{tame}}\widehat{\ot}_\Op \Op\ps{\Ga}$, where $\Ga= \Gal(\Q_\cyc/\Q)$. From this, we obtain an identification
\[T_{\mathrm{HL}} \cong T_\mathcal{F}\widehat{\ot}_\Op T_{\mathrm{tame}}^* \widehat{\ot}_\Op \Op\ps{\Ga}^{\iota}  \cong T_\mathcal{F}\widehat{\ot}_\Op T_{\mathrm{tame}}^* \widehat{\ot}_{\widetilde{R}} \widetilde{R}\ps{\Ga}^{\iota},\]
where here $T_{\mathrm{tame}}^*:= \Hom_{R_{\mathrm{tame}}}(T_{\mathrm{tame}},R_{\mathrm{tame}}(1))$ and $\widetilde{R}:=R_\mathcal{F}\ot R_{\mathrm{tame}}$. Note that $\widetilde{R}$ in this context is a power series ring in 3 variables. By an argument similar to that in Lemma \ref{Selcomplex cyc}, there is a natural isomorphism
\[H^i(SC(T_\mathcal{F}\widehat{\ot}_\Op T_{\mathrm{tame}}^*/\Q_\cyc))\cong H^i(SC(T_{\mathrm{HL}}/\Q))\]
for all $i$.

%If we write $A_{\mathrm{HL}}$ and $\widetilde{A}_{\mathrm{HL}}$ for the discrete modules of $T_{\mathrm{HL}}$ and $T_\mathcal{F}\widehat{\ot}_\Op T_{\mathrm{tame}}^*$ respectively, a Shapiro lemma argument similar to that in Lemma \ref{Sel Shapiro} will show that $\Sel_{str}(A_{\mathrm{HL}}/\Q) \cong \Sel_{str}(\widetilde{A}_{\mathrm{HL}}/\Q_\cyc)$.

We then consider the datum $(T_\mathcal{F}\ot T_{\mathrm{tame}}^*, T_{\mathcal{F},p}\ot T_{\mathrm{tame}}^*)_{R,\Q}$, and we leave the straightforward verification of \textbf{(C1)-(C4)} to the readers.  Suppose that $\p_f$ is a prime ideal of $R_{\mathcal{F}}$ of height 1 and that $\p_g$ is a prime ideal of $R_{\mathrm{univ}}$ of height 2 such that $T_f:=T_\mathcal{F}/\p_f$ and $T_g: = T_{\mathrm{univ}}/\p_g$ correspond to the Galois representations of cuspidal modular forms $f$  and $g$ with weights $k$ and $l$ respectively, where $k>l\geq 2$. Via the identification $R_{\mathrm{univ}} \cong R_{\mathrm{tame}}\widehat{\ot}_\Op \Op\ps{\Ga}$, we have a prime ideal $\overline{\p}_g$ of $R_{\mathrm{tame}}$ such that $T_{\mathrm{tame}}/\overline{\p}_g \cong T_g\ot \Op(\kappa^t)$. Together with $\p_f$, we obtain a prime ideal $\widetilde{\p}_{f,g}$ of $\widetilde{R}$ such that $T_\mathcal{F}\ot T_{\mathrm{tame}}^*/\widetilde{\p}_{f,g} \cong T_f\ot T_g(t)$. We can now verify that  $\textbf{(R1)}$ and  $\textbf{(R2)}$ are satisfied by a similar specialization argument as before.

\section{Appendix} \label{appendix}

Let $f,g,h$ be normalized cuspidal Hecke eigenforms with respective weights $k_f, k_g, k_h\geq 2$. We shall choose a large enough finite extension $E/\Q$ so that the respective associated representations $V_f, V_g, V_h$ can be viewed as representations over $E$ (after base-changing). The goal of this appendix is to give some sufficient conditions to ensure that the tensor products $V_f\ot_E V_g$ and $V_f\ot_E V_g\ot_E V_h$ are irreducible $\Gal(\overline{\Q}/\Q)$-representations. This is required for our discussion in the previous section for the verification of \textbf{(R1)}.

We start with the case of two modular forms.

\bp \label{product modular forms}
Suppose that $f$ has no complex multiplication. In the event that $g$ also has no complex multiplication, assume further that $k_f\neq k_g$. Then $V_f\ot_E V_g$ is irreducible as a $\Gal(\overline{\Q}/\Q)$-representation.
\ep

\bpf
We first consider the case that $g$ has complex multiplication given by an imaginary quadratic field $K$. Let $c$ be an element in $\Gal(\overline{\Q}/\Q)$ that projects to the nontrivial element of $\Gal(K/\Q)$. By complex multiplication, there is an algebraic Hecke character $\Psi$ of $K$ of infinity type $(k_g-1,0)$ such that $V_g\cong \mathrm{Ind}^{\Gal(\overline{\Q}/\Q)}_{\Gal(\overline{\Q}/K)} E(\psi)$, where $\psi:\Gal(\overline{\Q}/K) \lra E^\times$ is the $p$-adic Hecke character associated $\Psi$. From this, we obtain
\[ V_f\ot_EV_g \cong  \mathrm{Ind}^{\Gal(\overline{\Q}/\Q)}_{\Gal(\overline{\Q}/K)}
\mathrm{Res}^{\Gal(\overline{\Q}/\Q)}_{\Gal(\overline{\Q}/K)}V_f(\psi).\]
Writing $W:= \mathrm{Res}^{\Gal(\overline{\Q}/\Q)}_{\Gal(\overline{\Q}/K)}V_f$, it then follows from Clifford theory that $V_f\ot_EV_g$ is reducible if and only if $W(\psi) \cong W(\psi)^c$. Now, as $W\cong W^c$, the latter would imply that $\psi^c\psi^{-1}$ is a quadratic character, and in particular, a finite order character. However, this cannot happen in view that $\Phi^c\Phi^{-1}$ has infinity type $(1-k_g, k_g-1)$. Thus, the representation $V_f\ot_E V_g$ must be irreducible.

We now come to the situation when $g$ has no complex multiplication. For this part of the discussion, we shall freely use the material and notation as given in the work of Loeffler \cite{Loe}. Let $H$ be the subgroup of $\Gal(\overline{\Q}/\Q)$ cut out by all the inner twists of $f$ and $g$ (see \cite[Definition 2.2.1]{Loe} for the definition of inner twists). By comparing determinants and taking into account that $k_f>k_g$, one can see that $f$ is not Galois-conjugate to a twist of $g$. Therefore, it follows from the discussion in \cite[Subsection 3.3]{Loe} (also see \cite[Theorem 3.4.2]{Loe}) that the image of $H$ under the group homomorphism
\[\rho_f\times \rho_g : \Gal(\overline{\Q}/\Q)\lra \mathrm{GL}(V_f)\times \mathrm{GL}(V_g)\]
is an open subgroup of
\[ \mathcal{G}:=\left\{(x_f,x_g,\la)\in B_f^\times \times B_g^\times \times \Q_p^\times ~:~ \mathrm{norm}(x_f) = \la^{1-k_f},~\mathrm{norm}(x_g) = \la^{1-k_g}\right\},  \]
where here $B_f$ and $B_g$ are the quaternion algebras defined over an appropriate localization of the respective inner twist fields of $f$ and $g$. (In particular, we note that the quaternion algebras presented here are not those in \cite[Definition]{Loe} but rather appropriate localizations of those.) Note that the group homomorphism
\[ \rho_{f\ot g}:\Gal(\overline{\Q}/\Q)\lra \mathrm{GL}(V_f\ot_E V_g) \]
is precisely the composite of $\rho_f\times \rho_g$ with the Kronecker product $\mathrm{GL}(V_f)\times \mathrm{GL}(V_g)\lra \mathrm{GL}(V_f\ot_E V_g)$. Therefore, it suffices to show that $V_f\ot_E V_g$ is irreducible when viewed as a $\mathcal{G}$-representation which in turn reduces to verifying that $V_f\ot_E V_g$ is irreducible when viewed as a $\mathrm{Lie}(\mathcal{G})$-representation, where $\mathrm{Lie}(\mathcal{G})$ denotes $\mathcal{G}$'s Lie algebra. Note that $\mathcal{G} = \mathfrak{sl}_1(B_f)\oplus \mathfrak{sl}_1(B_g)\oplus \Qp$, where the final $\Qp$-summand acts on $V_f\ot_E V_g$ via scalar multiplication. Thus, it remains to show that $V_f\ot_EV_g$ is irreducible as a $\mathfrak{h}$-representation, where $\mathfrak{h}= \mathfrak{sl}_1(B_f)\oplus \mathfrak{sl}_1(B_g)$.

Enlarging $E$, if necessary, we may assume that the quaternion algebras $B_f$ and $B_g$ split. Consequently, we have $\mathfrak{h}\cong \mathfrak{sl}_2(E)\oplus \mathfrak{sl}_2(E)$ and $V_f\ot_EV_g\cong \mathrm{Std}_2(E)\boxtimes \mathrm{Std}_2(E)$, where $\mathrm{Std}_2$ is the standard representation of $\mathfrak{sl}_2$.
We then have the following calculations:
\begin{eqnarray*}
  \mathrm{End}_{\mathfrak{h}}\big(\mathrm{Std}_2(E)\boxtimes \mathrm{Std}_2(E)\big) \!&\cong &\! \mathrm{End}_{\mathfrak{sl}_2(E)}\big(\mathrm{Std}_2(E))\ot_E \mathrm{End}_{\mathfrak{sl}_2(E)}\big(\mathrm{Std}_2(E)\big) \\
   &=&\! E\ot_EE = E.
\end{eqnarray*}
From which, we see that $V_f\ot_E V_g$ is irreducible as a $\mathfrak{h}$-representation. This proves what we set to show.
\epf

For three modular forms, we have the following.

\bp \label{product three modular forms}
Suppose that $f,g,h$ have no complex multiplication and that their weights are mutually distinct. Then $V_f\ot_E V_g\ot_E V_h$ is irreducible as a $\Gal(\overline{\Q}/\Q)$-representation.
\ep

\bpf
By the weights assumption, a determinant comparison will show that $f,g,h$ are mutually not Galois-conjugate to one another. Therefore, we may apply \cite[Theorem 3.4.2]{Loe} to see that the image of a certain finite-index subgroup of $\Gal(\overline{\Q}/\Q)$ under
the group homomorphism
\[\rho_f\times \rho_g \times \rho_h : \Gal(\overline{\Q}/\Q)\lra \mathrm{GL}(V_f)\times \mathrm{GL}(V_g) \times \mathrm{GL}(V_h)\]
is an open subgroup of
\[\left\{(x_f,x_g,x_h,\la)\in B_f^\times \times B_g^\times \times B_h^\times \times \Q_p^\times ~:~ \mathrm{norm}(x_f) = \la^{1-k_f},~\mathrm{norm}(x_g) = \la^{1-k_g},~\mathrm{norm}(x_h) = \la^{1-k_h} \right\},  \]
where here $B_f, B_g$ and $B_h$ are certain appropriate quaternion algebras. Then we may proceed by a similar argument to that in the Proposition \ref{product modular forms}.
\epf

\end{document}